\documentclass[12pt, reqno, twoside]{amsart}
\usepackage{amsmath, amsthm, amscd, amsfonts, amssymb, graphicx, color}
\usepackage[bookmarksnumbered, colorlinks, plainpages]{hyperref}
\hypersetup{colorlinks=true,linkcolor=red, anchorcolor=green, citecolor=cyan, urlcolor=red, filecolor=magenta, pdftoolbar=true}
\usepackage{enumerate,amsfonts,mathtools,esint}
\usepackage{mathrsfs}
\usepackage{fancyhdr}
\usepackage{array}
\usepackage{hyperref}

\usepackage[toc,page]{appendix}
\usepackage{tabularx}
\usepackage{graphicx}

\usepackage{tikz}

\DeclareFontFamily{OT1}{rsfs}{}
\DeclareFontShape{OT1}{rsfs}{m}{n}{ <-7> rsfs5 <7-10> rsfs7 <10-> rsfs10}{}
\DeclareMathAlphabet{\mycal}{OT1}{rsfs}{m}{n}

\def\mcG{{\mycal G}}

\def\mcP{{\mycal P}}

\def\bA{{\mathbf{A}}}
\def\SO{{\bf SO}}

\def\bJ{{\mathbf{J}}}

\def\b0{{\mathbf{0}}}

\def\bH{\mathbf{H}}

\newcommand{\Q}{\bold{Q}}

\newcommand{\dive}{\textnormal{div}}
\newcommand{\la}{\langle}
\newcommand{\ra}{\rangle}
\newcommand{\In}{\bold{I}_n}
\newcommand{\msL}{\mathscr{L}}
\newcommand{\matrixexp}{\textnormal{exp}}

\newcommand{\R}{\mathbb{R}}
\newcommand{\Xn}{\mathbb{X}_n}

\newcommand{\An}{\mathbb{A}_n}

\newcommand{\mcH}{\mathscr{H}}

\newcommand{\sumell}{\sum_{\ell = 1}^N}
\newcommand{\dG}{\dot{\mathscr{G}}}
\newcommand{\dH}{\dot{\mathscr{H}}}

\newcommand{\m}{\mathsf{m}}
\newcommand{\wc}{\frac{2m_\ell\pi}{\mathsf{H}(b)}}
\newcommand{\rad}{\mathscr{R}_u}
\newcommand{\sph}{\mathscr{S}_u}
\newcommand{\A}{\mathsf{A}}
\newcommand{\B}{\mathsf{B}}

\newcommand{\uab}{[u;\mathsf{A},\mathsf{B}]}

\numberwithin{equation}{section}

\usepackage{enumitem}
\newtheorem{theorem}{Theorem}[section]
\newtheorem{proposition}{Proposition}[section]
\newtheorem{lemma}[theorem]{Lemma}
\newtheorem{corollary}[theorem]{Corollary}
\theoremstyle{definition}
\newtheorem{remark}[theorem]{Remark}
\numberwithin{equation}{section}

\begin{document}

\setcounter{page}{1}

\title[Nonlinear Elliptic System in Divergence Form]{On Multiple Solutions to a Family of Nonlinear Elliptic Systems in Divergence Form Coupled with an Incompressibility Constraint}

\author[A. Taheri]{Ali Taheri}
\author[V. Vahidifar]{Vahideh Vahidifar}

\address{School of Mathematical and  Physical Sciences, 
University of Sussex, Falmer, Brighton, United Kingdom.}
\email{\textcolor[rgb]{0.00,0.00,0.84}{a.taheri@sussex.ac.uk}} 
 
\address{Department of Applied Mathematics and Department of Mechanical Engineering, Iran University of Science and Technology, IUST, Tehran, Iran}
\email{\textcolor[rgb]{0.00,0.00,0.84}{engvahidifar@gmail.com}}



\begin{abstract}
The aim of this paper is to prove the existence of multiple solutions for a family of nonlinear elliptic systems in divergence form coupled 
with a pointwise gradient constraint:
\begin{align*}
\left\{
\begin{array}{ll}
\dive\{\A(|x|,|u|^2,|\nabla u|^2) \nabla u\} + \B(|x|,|u|^2,|\nabla u|^2) u 
= \dive \{ \mcP(x) [{\rm cof}\,\nabla u] \} \quad &\text{ in} \ \Omega , \\
\text{det}\, \nabla u = 1 \ &\text{ in} \ \Omega , \\
u =\varphi \ &\text{ on} \  \partial \Omega, 
\end{array}
\right. 
\end{align*}
where $\Omega \subset \mathbb{R}^n$ ($n \ge 2$) is a bounded domain, $u=(u_1, \dots, u_n)$ is a vector-map and $\varphi$ is a 
prescribed boundary condition. Moreover $\mathscr{P}$ is a hydrostatic pressure associated with the constraint $\det \nabla u \equiv 1$ 
and $\A = \A(|x|,|u|^2,|\nabla u|^2)$, $\B = \B(|x|,|u|^2,|\nabla u|^2)$ are sufficiently regular scalar-valued functions satisfying suitable 
growths at infinity. The system arises in diverse areas, e.g., in continuum mechanics and nonlinear elasticity, as well as geometric function 
theory to name a few and a clear understanding of the form and structure of the solutions set is of great significance. The geometric type of solutions 
constructed here draws upon intimate links with the Lie group ${\bf SO}(n)$, its Lie exponential and the multi-dimensional curl operator 
acting on certain vector fields. Most notably a discriminant type quantity $\Delta=\Delta(\A,\B)$, prompting from the PDE, will be shown 
to have a decisive role on the structure and multiplicity of these solutions.
\end{abstract}
 \maketitle


\section{Introduction}

This paper is motivated by questions on the existence and multiplicity of solutions to the following family of nonlinear elliptic 
systems in divergence form coupled with a pointwise gradient (incompressibility) constraint: 
\begin{align} \label{PDE1} 
\left\{
\begin{array}{ll}
\dive {\bf A} (x, u, \nabla u) + {\bf B} (x,u,\nabla u)  = \dive \{\mcP(x) [{\rm cof}\,\nabla u]\} 
\ &\text{in} \ \Omega,\\
\text{det} \nabla u = 1 \ &\text{in} \ \Omega, \\
u =\varphi \ &\text{on} \  \partial \Omega. 
\end{array}
\right. 
\end{align}
Here $\Omega \subset \mathbb{R}^n$ (with $n \ge 2$) is a bounded domain having a sufficiently smooth boundary $\partial \Omega$, 
${\bf A} = {\bf A}(x,u,\nabla u)$ and ${\bf B} = {\bf B}(x,u,\nabla u)$ are sufficiently regular $n \times n$ and $n \times 1$ matrix 
fields respectively, $u=(u_1, \dots, u_n)$ is an unknown vector-map defined on $\Omega$ with $\nabla u = [\partial u_i/\partial x_i : 1 \le i, j \le n]$ 
its gradient field, that is required to satisfy the pointwise incompressibility constraint $\det \nabla u \equiv 1$, and ${\rm cof}\,\nabla u$ denotes the 
cofactor matrix of $\nabla u$. In \eqref{PDE1} $\mcP$ is an {\it a priori} unknown scalar function (technically the Lagrange multiplier but also known 
as the hydrostatic pressure). Furthermore $\varphi$ is a prescribed boundary condition and the divergence operator acts on the rows of the 
two matrix fields on the left and right respectively.

The presence of the term $\dive \{\mcP [{\rm cof}\,\nabla u]\}$ is linked to the incompressibility constraint which makes the problem much harder 
compared to the unconstrained case where one typically either has no condition on the Jacobian determinant $\det \nabla u$ or an orientation 
preserving (pointwise positivity) condition $\det \nabla u>0$. Note that in the unconstrained case the expression on 
the right in the first line is zero (or equivalently the hydrostatic pressure $\mcP$ is constant).

This system arises in various fields ranging from nonlinear elasticity and continuum mechanics to geometric function theory 
({\it see} \cite{Ant, Ball1, CP, GiaModSou, HLeD} and the references therein) and the fundamental 
problem here is to establish the existence of solution pairs $(u, \mcP)$ subject to prescribed Dirichlet (or what is often called 
in the elasticity context as pure displacement) boundary conditions $u \equiv \varphi$ on $\partial \Omega$. For background and 
motivation {\it see} \cite{Ant, Ball1, Ball4, CP, MC, RY, Taheri-books, TeR, Valent} and for further studies and works in this and closely 
related directions 
{\it see}  \cite{AIMO, BallSch, BOP, BevDean, CK, EG, FosR, IO, HLeD, LeTallec, ALICM1, ALICM2, ALICM3, MsonTa, MsonTb, ST}.

For the sake of clarity, by a solution to the system \eqref{PDE1} in this paper, 
we mean a pair $(u, \mathscr{P})$ where the vector-map $u=(u_1, \dots, u_n)$ is of class 
$\mathscr{C}^2(\Omega, \mathbb{R}^n) \cap \mathscr{C}(\overline\Omega, \mathbb{R}^n)$, $\mathscr{P}$ is of class 
$\mathscr{C}^1(\Omega) \cap \mathscr{C}(\overline \Omega)$ and the pair satisfy the system \eqref{PDE1} in the 
pointwise (classical) sense. If the choice of $\mathscr{P}$ is clear from the context, or when no explicit reference to it is needed, 
we often abbreviate by saying that $u$ is a solution.

In nonlinear elasticity where a form of \eqref{PDE1} is encountered, the system represents the equilibrium equations 
of an incompressible material occupying the region $\Omega \subset \mathbb{R}^n$ ($n=2$ or $3$) in its reference 
(undeformed) configuration. Solutions here are incompressible deformations for which the body under the action of 
the external body forces and the prescribed displacement boundary conditions is at {\it rest} (or equilibrium) and are 
thus of great physical significance. Note additionally that in the hyperelastic case, the system represent the Euler-Lagrange 
equation associated with the constrained total elastic energy [{\it see} the integral \eqref{I-W-energy} below] and solutions in this context 
are equilibria as well as energy extremisers. (See below for more and Section \ref{Derivation-EL}. 
See also \cite{Ball1, BallSch, CK, CP, FosR, HLeD, LeTallec}.)

Now in order to further motivate \eqref{PDE1} and discuss the above in more detail consider a twice continuously differentiable 
stored energy density $W=W(x, u, \zeta)$ with $x \in \Omega$, $u \in \mathbb{R}^n$ and $\zeta \in \mathbb{R}^{n \times n}$ 
[or $\zeta$ in some fixed neighbourhood of the linear group ${\bf SL}(n)$]. 
For any incompressible deformation $u$ of $\Omega$, {\it i.e.}, any weakly differentiable map $u=(u_1, \dots, u_n)$ satisfying 
$\det \nabla u \equiv 1$ {\it a.e.} in $\Omega$, let its total elastic energy be given by the integral 
\begin{equation} \label{I-W-energy}
\mathbb{E}[u] = \int_\Omega W(x, u(x), \nabla u(x)) \, dx.  
\end{equation}
Incorporating boundary conditions and the growth of $W$ prompts one to introduce the space 
$\mathscr{A}^p_\varphi (\Omega) = \{ u \in \mathscr{W}^{1,p}(\Omega, \R^n) | \mbox{$\det \nabla u = 1$ {\it a.e.} 
in $\Omega$, $u = \varphi$ on $\partial \Omega$} \}$ for suitable choice of $1\le p<\infty$. Here boundary values 
are interpreted in the sense of traces of Sobolev functions. We hereafter refer to 
$\mathscr{A}^p_\varphi (\Omega)$ as the space of admissible incompressible maps or deformations.

The Euler-Lagrange equation associated with the energy integral \eqref{I-W-energy} over the space of admissible 
maps $\mathscr{A}^p_\varphi (\Omega)$ can be formally derived (see Section \ref{Derivation-EL} at the end) 
and seen to be given by the nonlinear system,  
\begin{equation} \label{PDE1-var}
-\dive \{ W_\zeta (x, u, \nabla u) - \mcP [{\rm cof}\,\nabla u] \} + W_u(x,u,\nabla u) =0, 
\end{equation} 
where $W_\zeta = [\partial W/\partial \zeta_{ij} : 1 \le i,j \le n]$ and $W_u = [\partial W/ \partial u_i : 1 \le i \le n]$. 
This system is evidently in the form \eqref{PDE1} with ${\bf A} = W_\zeta$ and ${\bf B}=-W_u$. 
Note however that \eqref{PDE1} is more general than \eqref{PDE1-var} in that there 
need not be any inherent relations between ${\bf A}$ and ${\bf B}$ in \eqref{PDE1}
whereas in the variational case leading to \eqref{PDE1-var} we have ${\bf A}_u=-{\bf B}_\zeta$ 
(specifically, $\partial {\bf A}_{ij}/\partial u_k = -\partial {\bf B}_k/\partial \zeta_{ij}$ with $1 \le i,j,k \le n$). 
In passing let us also note that by using the Piola identity ({\it see}, e.g., \cite{Ball1, CP}) and recalling the 
assumed regularity of solution pairs $(u, \mcP)$ we can write \eqref{PDE1-var} as 
\begin{equation} \label{PDE1-var-with-Piola}
-\dive \{ W_\zeta (x, u, \nabla u) \} + W_u(x,u,\nabla u) + [{\rm cof}\,\nabla u] \nabla \mcP = 0. 
\end{equation}

The system \eqref{PDE1-var} can be independently derived using the Lagrange multiplier method in the context of infinite 
dimensional differentiable manifold of incompressible maps ({\it cf.} \cite{HLeD} for details). An easy inspection here then 
shows that \eqref{PDE1-var} is also the Euler-Lagrange equation associated with the unconstrained energy integral 
$\mathbb{E}_\mathscr{P}$ incorporating the Lagrange multiplier and the constraint (which we leave the formal 
verification to the reader) given by 
\begin{equation} \label{I-P-W-energy}
\mathbb{E}_\mathscr{P} [u] = \int_\Omega \{ W(x, u(x), \nabla u (x)) - \mathscr{P}(x) [\det \nabla u (x) -1] \} \, dx. 
\end{equation}
Here evidently for any $u$ in $\mathscr{A}^p_\varphi (\Omega)$ we have $\mathbb{E}_\mathscr{P}[u]=\mathbb{E}[u]$. 
Let us point out that due to the {\it a priori} unknown regularity of the pressure field $\mathscr{P}$, and integrability of the 
Jacobian determinant $\det \nabla u$ the unconstrained energy integral $\mathbb{E}_\mathscr{P}$ in \eqref{I-P-W-energy} 
need not be everywhere well-defined, let alone, continuously Frechet differentiable on Sobolev spaces 
$\mathscr{W}^{1,p}(\Omega, \mathbb{R}^n)$. As a result standard tools from critical 
point theory do not carry over immediately to this setting ({\it cf.} \cite{P1, P2, PSM}) and so for the construction of energy extremisers (or critical 
points) other approaches and ideas are needed ({\it cf.} \cite{KT, TA1}). 

For the sake of this paper we focus on the case where the nonlinearities take the forms 
${\bf A} = \mathsf{A}(|x|,|u|^2,|\nabla u|^2)\nabla u$ and ${\bf B} = \B(|x|,|u|^2,|\nabla u|^2) u$ with 
$\A=\A(r,s,\xi)$ and $\B=\B(r,s,\xi)$ being sufficiently regular scalar-valued functions. This is called the {\it isotropic} case. [Note that in the setting 
of \eqref{I-W-energy} the latter amount to $W_\zeta = \mathsf{A}(|x|,|u|^2,|\nabla u|^2)\nabla u$ and 
$W_u=-\B(|x|,|u|^2,|\nabla u|^2) u$ where by writing $W(x,u,\nabla u)=F(|x|,|u|^2, |\nabla u|^2)$, it follows that 
$\A=\mathsf{F}_\xi$ and $\B=-\mathsf{F}_s$ (and so $\A_s+\B_\xi \equiv 0$). In the general case however there are no assumptions 
or relations linking $\A, \B$ and apart from standard regularity and growth (see below) the choices of $\A$ 
and $\B$ are independent and arbitrary.]

Now in view of the structure assumptions on the nonlinearities in place, the assumed regularity of solution pairs $(u, \mcP)$, 
and an application of Piola identity, the system in \eqref{PDE1} can be re-written in the form, 
\begin{align}\label{pre-EL1}
\left\{
\begin{array}{ll}
\dive[\A(|x|, |u|^2, |\nabla u|^2) \nabla u] + \B(|x|, |u|^2, |\nabla u|^2) u  = [{\rm cof}\,\nabla u] \nabla \mcP&\text{in} \ \Omega, \\
\text{det} \nabla u = 1 \ &\text{in} \ \Omega, \\
u =\varphi \ &\text{on} \  \partial \Omega, 
\end{array}
\right. 
\end{align}
Since $\det \nabla u \equiv 1$ we have $\det {\rm cof} \, \nabla u \equiv 1$ and 
$[{\rm cof}\,\nabla u]^{-1} = [\nabla u]^t$ 
and so we can write the {\it constrained} system \eqref{pre-EL1} in the more tractable gradient form 
\begin{align}\label{EL1}
\Sigma[(u, \mathscr{P}); (\mathsf{A}, \mathsf{B})] = 
\left\{
\begin{array}{ll}
[\nabla u]^t \left\{ \dive[\A\nabla u] +\B u\right\} = \nabla \mcP&\text{in} \ \Omega, \\
\text{det} \nabla u = 1 \ &\text{in} \ \Omega, \\
u =\varphi \ &\text{on} \  \partial \Omega. 
\end{array}
\right. 
\end{align}
where we have written $\A=\A(|x|, |u|^2, |\nabla u|^2)$ and $\B=\B(|x|, |u|^2, |\nabla u|^2)$ for brevity. 
It is also convenient to abbreviate the PDE in the first line of \eqref{EL1} to $\msL\uab = \nabla\mcP$ 
by introducing the differential operator 
\begin{align} \label{Lintro}
\msL\uab &= [\nabla u]^t\left\{\dive [\A\nabla u] + \B u \right\} \nonumber \\
&= [\nabla u]^t[\nabla u]\nabla \A + \A[\nabla u]^t\Delta u + \B[\nabla u]^t u.
\end{align}
Our aim here is to establish the existence of multiple solutions to the nonlinear system \eqref{EL1}. We confine 
to the geometric setting where the domain is an $n$-annulus, for definiteness, $\Omega=\{a<|x|<b\}$ with $b>a>0$ 
and $\varphi$ is the identity map. In this context a class of incompressible smooth maps with geometric significance 
are introduced and shown to lead to an infinitude of solutions. For related works on non-uniqueness in the 
incompressible setting {\it see} \cite{Bev, BevDean, MTa, ALICM1, ALICM2, MsonTa, MsonTb} and for results 
on uniqueness {\it see} \cite{KnS, ShTconv} (see also \cite{Ant, CP, HLeD, Valent}).

The study of multiple solutions 
to the constrained system $\Sigma[(u, \mcP); \A, \B]$, by way of construction, links here to a closely related 
unconstrained system for a vector-map $\mathsf{f}=(f_1, \dots, f_d)$, in a set of new variables $y=(y_1, \dots, y_N)$, 
and on a new region $\An \subset \mathbb{R}^N$, given by ({\it see} Section \ref{Restricted and Curl Section} for details)
\begin{align} 
\left\{
\begin{array}{ll}
\dive \left[ \mathscr{A}_i (y, \nabla \mathsf{f}) \nabla f_i \right]  = 0 \ &\text{ in } \An , \\
\mathsf{f} \equiv \mathsf{g} \ &\text{ on } (\partial\An)_\mathsf{D}, \qquad 1 \le i \le d, \\
\mathscr{A}_i(y, \nabla \mathsf{f}) \partial_\nu f_i= 0 \ &\text{ on } (\partial\An)_\mathsf{N}.
\end{array}
\right. 
\end{align}
Here $\mathsf{g} = (g_1,\dots,g_d)$ is a map defined on the so-called Dirichlet part of the boundary 
$(\partial\An)_\mathsf{D}$ ({\it see} below) describing the boundary values of the vector-map $\mathsf{f}$ itself whilst on the 
Neumann part $(\partial \An)_\mathsf{N} = \partial\An \setminus (\partial\An)_\mathsf{D}$ (the remainder of $\partial \An$) 
$\mathsf{f}$ is free. Additionally
\begin{equation}
\mathscr{A}_i (y, \nabla \mathsf{f}) = y_i^2 \A \bigg(z, z^2, n+\sum_{j=1}^d y_j^2 |\nabla f_j|^2 \bigg) \mathscr{J}(y), 
\qquad 1 \le i \le d, 
\end{equation}
with $\mathscr{J}(y) = y_1\cdots y_d$ and $z = \|y\|$ denoting the Euclidean $2$-norm of the $N$-vector $y$. 
The existence and multiplicity of solutions to this unconstrained system is discussed in Sections 
\ref{Restricted and Curl Section}-\ref{Non variational WD Sec} and the crucial connection between the two 
systems proved in Proposition \ref{connection systems} and its two corollaries Propositions 
\ref{whirl full vs restrcited} and \ref{reduction}. The main existence and multiplicity results of the paper 
are then presented in Theorem \ref{Dirsoln} and Theorem \ref{multiplicity even F}. As is apparent from 
the analysis in Section \ref{Non variational WD Sec}, a discriminant like object $\Delta=\Delta(u; \A, \B)$, 
plays a crucial role in the structure and dimensional parity of solutions. Let us end this introduction by 
formalising the assumptions on $\A, \B$ and fixing some key notation and terminology. 
\qquad \\
{\bf Assumptions on $\mathsf{A}, \mathsf{B}$.} We assume $\A = \A(r,s,\xi)$, $\B = \B(r,s,\xi)$ to be 
of class $\mathscr{C}^1(U)$, where $U = U[a, b] = [a, b] \times ]0, \infty[ \times ]0, \infty[$ 
with $\A > 0$, $\A_\xi \geq 0$ for all $(r, s, \xi) \in U$ and that for every compact set $K \subset ]0, \infty[$ there are constants 
$c_1 =c_1(K), \, c_2 =c_2(K)>0$ such that $c_1|\zeta|^{p-1} \leq \A(r, s, \zeta) |\zeta| \leq c_2|\zeta|^{p-1}$ for all 
$(r, s, \zeta) \in U, \, s \in K$ and $p > 1$. \\
{\bf Notation.} Throughout the paper we write $|x| = r$ and $\Theta = \nabla|x|= x|x|^{-1}$. By ${\bf I}={\bf I}_n$ we denote the $n\times n$ 
identity matrix. We write ${\bf J}=\sqrt{-{\bf I}_2}$ for the $2 \times 2$ skew-symmetric matrix with ${\bf J}_{12}=-1$ and write 
${\bf R}[\alpha] = {\rm exp} \{\alpha {\bf J}\}$ for the ${\bf SO}(2)$ matrix of rotation by angle $\alpha \in \mathbb{R}$ (in particular ${\bf J}={\bf R}[\pi/2]$).  
We write $y=(y_1, \dots, y_N)$ for the vector of $2$-plane radial variables associated with $x=(x_1, \dots, x_n)$ defined as follows: 
when $n = 2N$ we set $y_\ell = (x_{2\ell-1}^2 + x_{2\ell}^2)^{1/2}$ for $1 \leq \ell \leq N$ and when $n = 2N-1$ we set $y_\ell$ as 
before for $1 \leq \ell \leq N-1$ and $y_N = x_n$. For $b>a>0$ we write $\Xn=\{x \in \mathbb{R}^n : a<|x|<b\}$ and set 
$\An = \{ y \in \mathbb{R}^N_+ : a < \|y\| < b \}$ when $n=2N$ and $\An= \{ y \in \mathbb{R}^{N-1}_+ \times \R : a < \|y\| < b \}$ 
when $n=2N-1$. Here $|x|=(x_1^2 + \dots + x_n^2)^{1/2}$ and $\|y\| = (y_1^2 + \dots + y_N^2)^{1/2}$ denote the $2$-norms of 
the $n$-vector $x$ and $N$-vector $y$ respectively. In either case we have $\Xn \subset \mathbb{R}^n$ and $\An \subset \mathbb{R}^N$.  
Finally we write $(\partial\An)_\mathsf{D} = \{y \in \partial \mathbb{A}_n: \|y\|=a\} \cup \{y \in \partial \mathbb{A}_n: \|y\|=b\}$ and 
$(\partial\An)_\mathsf{N} = \partial \An \backslash (\partial\An)_\mathsf{D}$. Thus here $\partial \An= (\partial\An)_\mathsf{D} \cup (\partial\An)_\mathsf{N}$. 
It is often convenient to write $z=\|y\|$ and $\mathsf{1}=\mathsf{1}_d=(1, \dots, 1)$ for the $d$-vector whose 
components are all $1$. Vector inner product is denoted by $\langle u, v \rangle$ and matrix inner product by $E : F={\rm tr}(E^tF)$. 
Finally we use the standard notation for Sobolev spaces $\mathscr{W}^{1,p}$ (as, {\it e.g.}, \cite{Taheri-books}).


\section{The action $\msL\uab$ and the radial and spherical parts of $u$} \label{maps}
\setcounter{equation}{0}

Given a nowhere vanishing $u \in \mathscr{W}^{1,p}(\Omega, \R^n)$ (i.e., $u$ non-zero {\it a.e.} in $\Omega$) we decompose $u$ into a {\it radial} 
part $\rad$ and a {\it spherical} part $\sph$ by writing $\rad = |u|$ and $\sph = u|u|^{-1}$ respectively. A straightforward calculation then gives the 
gradients
\begin{equation}\label{gradRS}
\nabla \rad = \nabla |u| = \frac{[\nabla u]^tu}{|u|} , \qquad \nabla\sph = \nabla (u |u|^{-1}) = \left(\In - \frac{u}{|u|}\otimes \frac{u}{|u|}\right)\frac{\nabla u}{|u|} ,
\end{equation}
with $\In$ the $n\times n$ identity matrix. Of particular interest below are the two symmetric matrix fields relating to the left and right Cauchy-Green 
tensors (\cite{Ant, Ball1, CP}):
\begin{equation} \label{XY}
\mathbf{X}[u] = [\nabla u]^t[\nabla u] - \In , \qquad \mathbf{Y}[u] = [\nabla u][\nabla u]^t- \In. 
\end{equation}
Clearly these fields vanish {\it iff} $\nabla u$ is an orthogonal matrix {\it a.e.} in $\Omega$  
and so as such serve as a measure of the closeness of $\nabla u$ to the group ${\bf O}(n)$. Note also that 
$|\sph|=1$ and so $[\nabla \sph]^t \sph=0$ whilst $\rad^2 |\nabla \rad|^2 = \la {\bf Y}[u] u, u \ra + \rad^2$. Further conclusions are as below.

\begin{lemma}\label{key1}
Suppose $\rad$, $\sph$ are as in \eqref{gradRS} and $\mathbf{X}[u]$, $\mathbf{Y}[u]$ are as in \eqref{XY}. Then the following relations hold:
\begin{flalign*}
(i)& \ \nabla u = \rad \nabla\sph + \sph \otimes \nabla \rad , & \\
(ii)& \ [\nabla u]^tu = \rad (\rad [\nabla \sph]^t + \nabla \rad \otimes \sph) \sph = \rad\nabla \rad , & \\
(iii)& \ \mathbf{X}[u] = \rad^2[\nabla \sph]^t[\nabla \sph] + \nabla\rad \otimes \nabla\rad - \In , & \\
(iv)& \ \mathbf{Y}[u] = \rad^2[\nabla \sph][\nabla\sph]^t + \rad \nabla \sph\nabla \rad \otimes \sph + & \\
&\qquad \qquad \qquad + \rad\sph \otimes \nabla\sph \nabla \rad + |\nabla \rad|^2 \sph \otimes \sph - \In , & \\
(v)& \ |\nabla u|^2 = tr\{[\nabla u]^t[\nabla u]\} = tr\{[\nabla u][\nabla u]^t\} = \rad^2|\nabla \sph|^2 + |\nabla \rad|^2 , & \\
(vi)& \ \mathbf{X}[u]\nabla(|u|^2) = 2 \rad (\rad^2[\nabla \sph]^t[\nabla \sph]\nabla \rad + |\nabla\rad|^2 \nabla \rad - \nabla \rad).
\end{flalign*}
\end{lemma}

\begin{lemma}\label{key2}
If $u$ is second order differentiable then with $\rad$, $\sph$ as in \eqref{gradRS} we have: 
\begin{flalign*}
(i)& \ \Delta u = \rad \Delta \sph + 2\nabla \sph \nabla \rad + \Delta \rad \sph , & \\
(ii)& \ [\nabla u]^t\Delta u = \rad^2[\nabla\sph]^t\Delta\sph + \rad [2 [\nabla\sph]^t[\nabla\sph] 
+\langle \sph,\Delta\sph\rangle {\bf I}_n] \nabla\rad + \Delta\rad \nabla \rad , & \\
(iii)& \ \nabla(|\nabla u|^2) 
= 2 \rad |\nabla \sph|^2 \nabla \rad + 2 \rad^2 \nabla^2 \sph \nabla \sph + 2 \nabla^2 \rad \nabla \rad. 
\end{flalign*}
\end{lemma}

\begin{proof}
These are all consequences of direct differentiation and routine calculations. 
\end{proof}


\begin{lemma}\label{Lu}
Suppose $u$ is second order differentiable and $\msL[u;\A,\B]$ is as in  \eqref{Lintro}. Then the following relation holds: 
\begin{align}\label{Lcoef}
\msL[u;\A,\B] -& \nabla \A = \, \mathbf{X}[u]\nabla \A + \A[\nabla u]^t\Delta u + \B[\nabla u]^t u \\
=& \, \A_\xi \mathbf{X}[u]\nabla(|\nabla u|^2) + \A_s \mathbf{X}[u]\nabla(|u|^2) + \A_r \mathbf{X}[u]\Theta+ \A [\nabla u]^t\Delta u + \B[\nabla u]^t u. \nonumber 
\end{align}
Here $\A_r = \A_r(r,s,\xi), \, \A_s = \A_s(r,s,\xi)$ and $\A_\xi = \A_\xi(r,s,\xi)$ denote the respective partial derivatives of $\A$ 
whilst $\B=\B(r,s,\xi)$. All arguments are at $(r,s,\xi)=(|x|, |u|^2, |\nabla u|^2)$.
\end{lemma}

\begin{proof}
This follows from \eqref{Lintro} after substituting for $\bold{X}[u]$ from \eqref{XY} and then rearranging terms.
\end{proof}


\qquad \\
{\bf Generalities on maps with $\rad =|x|$, $\sph=\Q \Theta$}: The class of maps we are interested in here are those whose 
radial and spherical parts are $\rad (x) = |x|$, $\sph (x) = \Q(x) \Theta$ respectively. Here $\Q$ is an $\SO(n)$-valued matrix 
field whose dependence on the spatial variables $x=(x_1, \dots, x_n)$ is through the $2$-plane radial variables $y=(y_1, \dots, y_N)$ 
described earlier.  Thus with a slight abuse of notation we hereafter write and think of $\Q=\Q(y)$ with $y=y(x)$ 
({\it see} \cite{ALICM1, ALICM2}).

We next define the set of $2N$ orthogonal $n$-vectors: $w^i = (0,\dots,0,x_{2i-1},x_{2i},0,\dots,0)$, 
$[w^i]^\perp = (0,\dots,0,-x_{2i},x_{2i-1},0,\dots,0)$ for $1 \leq i \leq d$; when $n=2d$ is even this completes 
the picture but when $n=2d+1$ is odd we set $w^N = (0,\dots,0,x_n)$, $[w^N]^\perp = (0,\dots,0)$. Hence 
$x = w^1 + \dots + w^N$, $\langle w^i, w^j \rangle = 0$, $\langle [w^i]^\perp, [w^j]^\perp \rangle =0$ for 
$1 \le i \neq j \le N$ and $\langle w^i, [w^j]^\perp \rangle = 0$ for all $1 \le i, j \le N$. Furthermore in relation 
to the variables $y_1,\dots,y_N$ introduced earlier we have $y_\ell = |w^\ell| = |[w^\ell]^\perp|$ when 
$1 \leq \ell \leq d$ noting that when $n=2d+1$ we have $w^N = (0,\dots,0,y_N)$ and $|y_N|=|w^N|=|x_n|$.

\begin{lemma}\label{ykey} 
For the $2$-plane radial variables $y = (y_1,\dots,y_N)$ we have: $\nabla y_\ell = w^\ell/y_\ell$, 
$\la \nabla y_\ell , \nabla y_k \ra = \delta_{\ell k}$ and $\Delta y_\ell = 1/y_\ell$ except for 
$n=2d+1$ where $\Delta y_N=0$. Here $\nabla, \Delta$ are taken with respect to the $x=(x_1, \dots, x_n)$ variables.    
\end{lemma}

\begin{proof}
These follow by straightforward differentiation and considering the cases corresponding to even and odd $n$ separately.
\end{proof}

\begin{lemma}\label{wrs}
Assume $u = \Q(y_1, \dots, y_N) x$ with the matrix field $\Q$ being of class $\mathscr{C}^1$. 
Then $\nabla \rad = x/|x|$, $\Delta\rad = (n-1)/|x|$ and 
\begin{equation}
\nabla\sph =  \frac{1}{||y||} \Q (\In - \Theta \otimes \Theta) + \sum_{\ell = 1}^N \partial_\ell\Q  \Theta \otimes \nabla y_\ell, 
\qquad \nabla \sph \nabla \rad = \sum_{\ell = 1}^N \la \nabla y_\ell, \Theta \ra \partial_\ell\Q  \Theta.
\end{equation}
Moreover if $\Q$ is of class $\mathscr{C}^2$ then $u$ is second order differentiable and 
\begin{equation}
\Delta\sph = \frac{1-n}{||y||^2} \Q \Theta + \frac{1}{||y||} \sum_{\ell=1}^N \bigg\{ [ 2\partial_\ell \Q\nabla y_\ell - 2 \langle \nabla y_\ell , \Theta \rangle \partial_\ell\Q \Theta ] 
+ \partial^2_\ell\Q \Theta + \Delta y_\ell\partial_\ell\Q \Theta \bigg\}.
\end{equation}
\end{lemma}


We next prove further identities associated with such maps in line with Lemmas \ref{key1} and \ref{key2}. Note that  
$\partial_\ell$ stands for partial differentiation with respect to $y_\ell$ whilst $\nabla$ and $\Delta$ as applied to the variables 
$y=(y_1, \dots, y_N)$ are all with respect to $x=(x_1,\dots,x_n)$. 

\begin{lemma}\label{key1w}
Assume $u = \Q(y_1, \dots, y_N) x$ with the matrix field $\Q$ being of class $\mathscr{C}^1$. Then the following identities hold:
\begin{flalign*}
(i)& \ \nabla u = \Q + \sumell \partial_\ell\Q x \otimes \nabla y_\ell, & \\
(ii)& \ [\nabla u]^t u = \bigg[ \Q^t + \sumell \nabla y_\ell \otimes \partial_\ell \Q x \bigg]\Q x = x, & \\
(iii)& \ |\nabla u|^2 = tr\{[\nabla u]^t[\nabla u]\} = tr\{[\nabla u][\nabla u]^t\} 
= n + \sumell \left[2 \la \Q^t\partial_\ell\Q x, \nabla y_\ell \ra + |\partial_\ell\Q x|^2\right].
\end{flalign*}
\end{lemma}

\begin{lemma}\label{key1w XY}
Under the assumptions of the previous lemma on $u$ and with ${\bf X}[u], {\bf Y}[u]$ denoting the matrix fields in \eqref{XY} the following identities hold:
\begin{flalign*}
(i)& \ \mathbf{X}[u] = \sumell \left[ \Q^t\partial_\ell\Q x \otimes \nabla y_\ell + \nabla y_\ell \otimes \Q^t \partial_\ell\Q x\right] 
+ \sumell \sum_{k=1}^N \langle\partial_\ell\Q x , \partial_k\Q x \rangle \nabla y_\ell \otimes \nabla y_k, & \\
(ii)& \ \mathbf{Y}[u] = \sumell \left[ \Q\nabla y_\ell \otimes \partial_\ell\Q x +\partial_\ell\Q x \otimes \Q\nabla y_\ell 
+ \partial_\ell\Q x \otimes\partial_\ell\Q x \right], & \\
(iii)& \ \mathbf{X}[u]\nabla(|u|^2) = 2\bigg[ \sumell \langle \nabla y_\ell , x\rangle \Q^t\partial_\ell \Q x 
+ \sumell \sum_{k=1}^N \langle \partial_\ell \Q x , \partial_k \Q x \rangle \langle \nabla y_k , x \rangle \nabla y_\ell \bigg].
\end{flalign*}
\end{lemma}

\begin{proof}
The identities $(i)$ and $(ii)$ above follow from $(iii)$ and $(iv)$ in Lemma~\ref{key1} and the identities in Lemma~\ref{wrs}. 
To conclude $(ii)$ we use the relation $\langle \Q^t\partial_\ell\Q x,x\rangle = 0$ resulting from skew-symmetry. The third 
identity follows at once by noting $|u|^2=\rad^2 = r^2$.
\end{proof}

\begin{lemma}\label{key2wP1}
Assume $u =\Q(y_1, \dots, y_N) x$ with the matrix field $\Q$ being of class $\mathscr{C}^2$. Then $u$ is second order differentiable and the following identities hold:
\begin{align*}
(i)&\ \Delta u  = \sum_{\ell=1}^N\left[\partial^2_\ell\Q x+\Delta y_\ell\partial_\ell\Q x+ 2\partial_\ell\Q \nabla y_\ell\right], & \\
(ii)& \ [\nabla u]^t \Delta u = \sumell \left[ \Q^t\partial^2_\ell\Q x + \Delta y_\ell \Q^t\partial_\ell\Q x + 2\Q^t\partial_\ell\Q \nabla y_\ell \right] & \\
& \qquad \qquad + \sumell \sum_{k=1}^N \left[ \langle \partial_\ell\Q x, \partial^2_k\Q x\rangle + \Delta y_k \langle \partial_\ell\Q x, \partial_k\Q x \rangle 
+ 2\langle \partial_\ell\Q x,\partial_k\Q \nabla y_k \rangle\right] \nabla y_\ell, & \\ 
(iii)& \ \nabla(|\nabla u|^2) =  2 \sumell \left[ \partial_\ell \Q^t \Q \nabla y_\ell + \nabla^2 y_\ell \Q^t \partial_\ell \Q x 
+ \partial_\ell\Q^t \partial_\ell\Q x \right] \nonumber \\
& \qquad \qquad \qquad + \sumell \sum_{k=1}^N 2 \left[ \la (\partial_k \Q^t \partial_\ell \Q + \Q^t \partial_{\ell k} \Q) x, \nabla y_\ell \ra 
+ \la \partial^2_{\ell k}\Q x , \partial_\ell\Q x \ra \right] \nabla y_k. &
 \end{align*}
\end{lemma}

\begin{proof}
For $(i)$ we use the identities in Lemma~\ref{wrs} together with the description of the Laplacian given 
in identity $(i)$ in Lemma~\ref{key2}. We then obtain $(ii)$ by pre-multiplying this with $[\nabla u]^t$ using 
the description of $\nabla u$ given by $(i)$ in Lemma~\ref{key1w}. For $(iii)$ by invoking $(iii)$ in Lemma 
\ref{key1w} and expanding the gradient directly on each term we have  
\begin{align}
\nabla \la \Q^t\partial_\ell\Q x, \nabla y_\ell \ra =& \partial_\ell \Q^t \Q \nabla y_\ell + \nabla^2 y_\ell \Q^t \partial_\ell \Q x \nonumber \\
& + \sum_{k=1}^N [\nabla y_k \otimes \nabla y_\ell](\partial_k \Q^t \partial_\ell \Q + \Q^t \partial_{\ell k} \Q) x, 
\end{align}
and likewise  
\begin{align} \label{core term}
\nabla |\partial_\ell\Q x|^2 = \nabla \la \partial_\ell\Q x, \partial_\ell\Q x 
\ra = 2 \partial_\ell\Q^t\partial_\ell \Q x + \sum_{k=1}^N 2 [\nabla y_k \otimes \partial^2_{\ell k}\Q x] \partial_\ell \Q x.  
\end{align}
Putting these together and rearranging terms gives at once the desired conclusion. 
\end{proof}

\section{Whirls, maximal tori, and the block diagonal ${\bf SO}(n)$-valued matrix fields $\Q[\mathsf{f}]$}
\label{Whirls and Blocks Section}

Returning to the decomposition of $u$ into its radial and spherical parts, and prompted by symmetry considerations, we now 
specialise to the class of maps $u$ whose $\SO(n)$-valued matrix field $\Q$ in the spherical part $\sph$ takes values on a 
fixed maximal torus ${\bf T}$ of $\SO(n)$ ({\it cf.} \cite{ALICM2} for more on this). As any two maximal tori on a compact Lie 
group are conjugate to one another, for our purposes, and without loss of generality, we take the canonical maximal torus 
${\bf T}$ of $2 \times 2$ block diagonal matrices ${\bf T}=\{diag({\bf R}_1, \dots, {\bf R}_d)\}$ for $n=2d$ and 
${\bf T}=\{diag({\bf R}_1, \dots, {\bf R}_d, 1)\}$ for $n=2d+1$. Here ${\bf R}_j={\bf R}[\alpha_j] \in {\bf SO}(2)$ 
($1 \le j \le d$) with $(\alpha_1, \dots, \alpha_d) \in \mathbb{R}^d$ (see \cite{KnA, VNJ} for background and more on Lie 
groups and representations).

The implication of this is that we will express $\Q$ as a similar block diagonal matrix with each block described by a suitable 
angle of rotation function $f=f_\ell(y)$ ($1 \le \ell \le d$). Specifically, this leads to the explicit descriptions of the 
${\bf SO}(n)$-valued matrix fields
\begin{align}\label{Qeven}
\Q[\mathsf{f}](y) = \left( \begin{array}{ccccc}
{\bf R}[f_1(y)] & 0 & \dots & 0 & 0 \\
0 & {\bf R}[f_2(y)] & \dots & 0 & 0 \\
\vdots & \vdots & \ddots & \vdots & \vdots \\
0 & 0 & \dots & {\bf R}[f_{d-1}(y)] & 0 \\
0 & 0 & \dots & 0 & {\bf R}[f_d(y)]
\end{array}\right),
\end{align}
for when $n=2d$ and 
\begin{align}\label{Qodd}
\Q[\mathsf{f}](y) = \left( \begin{array}{cccccc}
{\bf R}[f_1(y)] & 0 & \dots  & 0 & 0 & 0\\
0 & {\bf R}[f_2(y)] & \dots  & 0 & 0 & 0\\
\vdots & \vdots & \ddots  & \vdots & \vdots & \vdots \\
0 & 0 & \dots &{\bf R}[f_{d-1}(y)] & 0 & 0 \\ 
0 & 0 & \dots & 0 & {\bf R}[f_d(y)] & 0 \\
0 &0 & \dots & 0 & 0 & 1
\end{array}\right),
\end{align}
for when $n=2d+1$. Hereafter we write $\mathsf{f}=\mathsf{f}(y)$ for the vector-map $\mathsf{f}=(f_1, \dots, f_d)$ and we 
refer to the resulting map $u=\Q[\mathsf{f}](y) x$ as a {\bf whirl map} or a {\bf whirl} for simplicity. It can be seen by 
direct verification that 
\begin{align} \label{QtQl}
\Q^t\partial_\ell \Q \nabla y_k = y_k^{-1} \partial_\ell f_k [w^k]^\perp, \quad 
\Q^t \partial_\ell \Q \, x = \sum_{i=1}^d \partial_\ell f_i [w^i]^\perp, 
\end{align}
and likewise 
\begin{align} \label{QtQl-next}
\partial_\ell \Q^t \partial_\ell \Q \, x = \sum_{i=1}^d (\partial_\ell f_i)^2w^i, \qquad 
\Q^t\partial^2_{\ell k} \Q \, x = \sum_{i=1}^d \left[ \partial^2_{\ell k} f_{i} [w^i]^\perp - \partial_\ell f_{i} \partial_k f_{i}  w^i\right], 
\end{align}
(with $1 \le \ell \le N$ and $1 \le k \le N$). Next using the first identity in \eqref{QtQl-next} we have 
\begin{align*} 
\langle \nabla y_\ell, \partial_k \Q^t \partial_k \Q \, x \rangle 
= \left\langle \nabla y_\ell, \sum_{i=1}^d (\partial_k f_i)^2w^i \right\rangle, 
\end{align*}
which then upon making note of the inner product relation $\la \nabla y_\ell , w^i \ra = y_\ell \delta_{i \ell}$ leads to 
\begin{align}
\la \nabla y_\ell , \partial_k \Q^t \partial_k \Q x \ra = 
\begin{dcases}
y_\ell (\partial_k f_\ell)^2 &\mbox{$1 \le \ell \le d$},\\
0 &\mbox{$\ell=N$, $n$ odd}.
\end{dcases}
\end{align}

\begin{lemma}\label{key2wP2}
Let $u$ be a whirl as defined above with matrix field $\Q$ of class $\mathscr{C}^2$. Then 
\begin{align} \label{core X term}
\mathbf{X}[u] \nabla(|\nabla u|^2) =& 2 \sum_{k=1}^N \sumell \left[ \langle \nabla y_k, \partial_\ell \Q^t\partial_\ell \Q x \rangle 
+ \langle \partial^2_{\ell k} \Q x , \partial_\ell \Q x \rangle\right] \times \nonumber \\
& \qquad \bigg[ \Q^t\partial_k \Q x +\sum_{i=1}^N  \langle \partial_k \Q x, \partial_i \Q x \rangle \nabla y_i \bigg].
\end{align}
\end{lemma}

\begin{proof}
Since $\la \Q^t\partial_\ell\Q x, \nabla y_\ell \ra=0$ we have $\nabla(|\nabla u|^2) = \sumell \nabla |\partial_\ell \Q x|^2$. 
Utilising \eqref{core term}, pre-multiplying by $\mathbf{X}[u]$ using $(i)$ in Lemma \ref{key1w XY} and \eqref{QtQl}, 
\eqref{QtQl-next} give the result. 
\end{proof}

We now turn to formulating the action of the differential operator $\msL[u;\A,\B]$ on a whirl $u$ whose radial and spherical parts have the forms 
$\rad(x)=|x|$, $\sph(x) = \Q[\mathsf{f}](y) \Theta$.

\begin{proposition}\label{LQprop}
Suppose $u$ is a whirl with matrix field $\Q$ of class $\mathscr{C}^2$. The action of the differential operator $\msL$ 
on $u$ can be reformulated as 
\begin{align}\label{LQ}
\msL\,[u;\A,\B] =& \nabla \A + \B x + 2\A_\xi \sum_{k=1}^N \sumell  \left[ \langle \nabla y_k, \partial_\ell \Q^t\partial_\ell \Q x \rangle 
+  \langle \partial^2_{\ell k} \Q x , \partial_\ell \Q x \rangle\right] \times \\
&\quad\times \bigg[ \Q^t\partial_k \Q x +\sum_{i=1}^N  \langle \partial_k \Q x, \partial_i \Q x \rangle \nabla y_i \bigg] \nonumber \\
&+ \left[2\A_s+ r^{-1}\A_r\right] \sumell\bigg\{ \langle \nabla y_\ell , x\rangle \Q^t\partial_\ell\Q x 
+ \sum_{k=1}^N \langle \partial_\ell\Q x, \partial_k \Q x\rangle \langle \nabla y_k ,x \rangle \nabla y_\ell \bigg\} \nonumber \\
&+ \A\sumell \bigg\{ \left[ \Q^t\partial^2_\ell\Q x + \Delta y_\ell \Q^t\partial_\ell\Q x + 2\Q^t\partial_\ell\Q \nabla y_\ell \right] \nonumber \\
&+\sum_{k=1}^N \left[ \langle \partial_\ell\Q x, \partial^2_k \Q x\rangle + \Delta y_k \langle\partial_\ell\Q x, \partial_k \Q x \rangle 
+ 2\langle \partial_\ell\Q x, \partial_k\Q \nabla y_k\rangle\right] \nabla y_\ell \bigg\} .  \nonumber 
\end{align}
The arguments of $\mathsf{A} = \A(r,s,\xi), \, \B = \B(r,s,\xi)$ and all subsequent derivatives in \eqref{LQ} are 
$(r,s,\xi) = (|x|,|u|^2,|\nabla u|^2) = (r, r^2 , n+ \sum_{\ell=1}^N |\partial_\ell\Q x|^2)$.
\end{proposition}

\begin{proof}
This follows by referring to \eqref{Lcoef}. Firstly, the coefficient of $\A_\xi$ is $\mathbf{X}[u]\nabla(|\nabla u|^2)$, 
given by identity \eqref{core X term} in Lemma~\ref{key2wP2}. The coefficient of $\A_s$ is $\bold{X}[u]\nabla(|u|^2)$, given by 
identity $(iii)$ in Lemma~\ref{key1w XY} and with $\nabla|x| = \Theta$, the coefficient of $\A_r$ is appropriately described 
above. Similarly, the coefficient of $\A$ is $[\nabla u]^t\Delta u$, described by identity $(ii)$ in Lemma~\ref{key2wP1} and 
by noting $(ii)$ in Lemma~\ref{key1w} we recover $\B x$ in \eqref{LQ}. 
\end{proof}

\begin{remark}{\em Using $(i)$ in Lemma $\ref{key1w}$ a whirl is seen to satisfy the incompressibility constraint. Indeed 
$\det\nabla u = \det [\Q + \sum_{\ell = 1}^N \partial_\ell\Q x \otimes \nabla y_\ell] = \det [\In + \sum_{\ell = 1}^N \Q^t\partial_\ell\Q x\otimes \nabla y_\ell]$. 
Now since for $p_i = \Q^t\partial_i\Q x$, $q_j = \nabla y_j$ we have $\langle p_i, q_j\rangle = 0$ for all $1 \leq i,j \leq N$ it follows 
that $\det[\In + \sum_{i=1}^N p_i \otimes q_i ] = 1$ $(${\it cf.} Lemma $3.1$ in \cite{ALICM1}$)$ and so $\det\nabla u = 1$ as claimed.}
\end{remark}

\begin{remark}{\em The boundary condition $u \equiv x$ on $\partial \Xn$ $($equivalently $\Q[\mathsf{f}] \equiv {\bf I}_n$ on $\partial \Xn$$)$ 
translates to $\mathsf{f} \equiv 0$ on $\{z=a\}$ and $\mathsf{f} \equiv 2 \mathsf{m} \pi$ at $\{z=b\}$ 
with $\mathsf{m}=(m_1, \dots, m_d) \in \mathbb{Z}^d$. This follows by observing that $x \in \partial \Xn \iff y=y(x) \in (\partial \An)_D$ 
with the segments $\{|x|=a\}$ and $\{|x|=b\}$ of $\partial \Xn$ corresponding to the segments $\{z=a\}$ and $\{z=b\}$ of $(\partial \An)_D$ 
respectively whilst ${\bf R}[\alpha] = {\bf I}_2 \iff \alpha=2m \pi$.}
\end{remark}


\section{An auxiliary system and the interrelation of two differential operators} 
\label{Restricted and Curl Section}

We begin the section by introducing a nonlinear unconstrained system in divergence form 
\begin{align} \label{ELwhirl}
\left\{
\begin{array}{ll}
\dive \left[ \mathscr{A}_i (y, \nabla \mathsf{f}) \nabla f_i \right]  = 0 \ &\text{ in } \An , \\
\mathsf{f} \equiv \mathsf{g} \ &\text{ on } (\partial\An)_\mathsf{D}, \qquad 1 \le i \le d, \\
\mathscr{A}_i (y, \nabla \mathsf{f}) \partial_\nu f_i = 0 \ &\text{ on } (\partial\An)_\mathsf{N} .
\end{array}
\right. 
\end{align}
Here $\mathsf{f} = (f_1,\dots,f_d)$ is the unknown vector with $\nabla \mathsf{f} = [\partial_k f_\ell : 1 \le \ell \le d, 1 \le k \le N]$, 
the divergence is taken with respect to the $(y_1, \dots, y_d)$ variables and the nonlinearity (coefficients) in the PDE are given by 
\begin{equation} \label{nonlinearity A def}
\mathscr{A}_i (y, \nabla \mathsf{f}) := \A \bigg(z, z^2, n+\sum_{\ell=1}^d y_\ell^2 |\nabla f_\ell|^2 \bigg) y_i^2 \mathscr{J}(y), 
\qquad z^2=\|y\|^2=\sum_{j=1}^N y_j^2, 
\end{equation}
where $\mathscr{J}(y) = y_1\cdots y_d$. Recall that 
$(\partial\An)_\mathsf{D} = \{y \in \partial \mathbb{A}_n: z=a\} \cup \{y \in \partial \mathbb{A}_n: z=b\}$, $\mathsf{g}=\mathsf{g}(y, \mathsf{m})$ 
is the piecewise constant map defined by $\mathsf{g}|_{z=a}\equiv 0$ and $\mathsf{g}|_{z=b} \equiv 2\m\pi$ with $\m \in \mathbb{Z}^d$ fixed 
whilst $(\partial\An)_\mathsf{N} = \partial \An \backslash (\partial\An)_\mathsf{D}$ and $\partial_\nu f_i = \nabla f_i \cdot \nu$ with $\nu$ the 
unit outward normal field on $(\partial\An)_\mathsf{N}$. The motivation for studying this system by way of its relation to $\msL\uab$ and the 
system \eqref{EL1} will become clear later on. First we establish the uniqueness of solutions to \eqref{ELwhirl}. 
We set $\mathscr{B}_p[\An;\mathsf{g}] 
= \{ \mathsf{f} = (f_1,\dots,f_d) \in \mathscr{W}^{1,p}(\An,\R^d) : \mathsf{f} \equiv \mathsf{g} \text{ on } (\partial\An)_\mathsf{D} \}$ 
and note that the unconstrained system \eqref{ELwhirl} here is {\it strictly} elliptic but not uniformly elliptic as a result of 
$\mathscr{J}(y) >0$ in $\An$ but $\mathscr{J}(y) \equiv 0$ on $(\partial \An)_\mathsf{N}$ [{\it see} \eqref{nonlinearity A def}]. 
Thus interestingly even the existence of solution falls outside standard theory.

\begin{proposition} \label{unique}
Given $\m = (m_1,\dots,m_d)\in \mathbb{Z}^d$ the solution $\mathsf{f}=(f_1, \dots, f_d) \in\mathscr{C}^2(\overline{\mathbb{A}}_n,\R^d)$ 
to the system \eqref{ELwhirl} is unique. \footnote{For some explicit examples of solutions to this system with the required degree of regularity 
see Sections \ref{Non variational WD Sec} and \ref{Infinitude in even d Sec}.}
\end{proposition}

\begin{proof} 
Let $\mathsf{f}^1, \mathsf{f}^2$ be two solutions to \eqref{ELwhirl} in $\mathscr{B}_p[\An;\m]$ 
and put $\mathsf{h}=\mathsf{f}^2-\mathsf{f}^1$. Then $\mathsf{h} \equiv 0$ on $(\partial\An)_\mathsf{D}$. 
Now using the monotonicity inequality $[\A(z,z^2,\zeta_2) - \A(z,z^2,\zeta_1)] (\zeta_2 - \zeta_1) \ge 0$ for 
$\zeta_1, \, \zeta_2 \in \R$ with $\zeta_1 = n+ \sum y_j^2 |\nabla f^1_j|^2, \, \zeta_2 = n+ \sum y_j^2 |\nabla f^2_j|^2$ 
(the sums over $1 \leq j \leq d$) it follows after multiplying by $\mathscr{J} \ge 0$ and substitution using \eqref{nonlinearity A def} that 
\begin{align}
0 \le &\sum_{\ell=1}^d y_\ell^2
\left[ \A \left(z,z^2,n + \sum_{j=1}^d y_j^2 |\nabla f^2_j|^2\right) - \A \left(z,z^2,n+ \sum_{j=1}^d y_j^2 |\nabla f^1_j|^2\right) \right] \times \nonumber \\
& \times \mathscr{J} \left( |\nabla f^2_\ell|^2 - |\nabla f^1_\ell|^2 \right) 
= \sum_{\ell=1}^d [\mathscr{A}_\ell (y, \nabla \mathsf{f}^2) - \mathscr{A}_\ell(y, \nabla \mathsf{f}^1)] (|\nabla f^2_\ell|^2 - |\nabla f^1_\ell|^2) \\
\le& \sum_{\ell=1}^d 2 \la \mathscr{A}_\ell (y, \nabla \mathsf{f}^2) \nabla f^2_\ell - \mathscr{A}_\ell(y, \nabla \mathsf{f}^1) \nabla f^1_\ell,  \nabla h_\ell \ra
- \sum_{\ell=1}^d [\mathscr{A}_\ell(y, \nabla \mathsf{f}^2) + \mathscr{A}_\ell(y, \nabla \mathsf{f}^1)] |\nabla h_\ell|^2. \nonumber 
\end{align}
Now integrating the above and taking advantage of $\mathsf{f}^1, \mathsf{f}^2$ being solutions to \eqref{ELwhirl} it follows after 
an application of the integration by parts formula and noting the vanishing of the integral of the expression on the second line above that 
\begin{equation}
\int_{\An} \sum_{\ell=1}^d - [\mathscr{A}_\ell(y, \nabla \mathsf{f}^1) + \mathscr{A}_\ell(y, \nabla \mathsf{f}^2)] |\nabla h_\ell|^2 \, dy \ge 0. 
\end{equation} 
As $\mathscr{A}_\ell>0$ inside $\An$ it follows by taking into account the connectedness of $\An$ and the Dirichlet boundary condition 
on $\mathsf{h}$ that $\mathsf{h} \equiv 0$. Thus $\mathsf{f}^1 \equiv \mathsf{f}^2$ as required. 
\end{proof}

We now aim to make the link between the unconstrained system \eqref{ELwhirl} 
and the PDE $\msL[u;\A,\B] = \nabla\mcP$ in the original system \eqref{EL1} more transparent. Towards this end, we begin by expanding 
the divergence in \eqref{ELwhirl} thus obtaining the formulation 
\begin{align}\label{divexp}
\frac{1}{\mathscr{J} y_i^2} & \dive[\mathscr{A}_i(y, \nabla \mathsf{f}) \nabla f_i] = \frac{1}{\mathscr{J} y_i^2} \sum_{k=1}^N \partial_k 
\bigg[\mathsf{A}\bigg(z,z^2,n+\sum_{j=1}^d y_j^2 |\nabla f_{j}|^2\bigg) y_i^2 y_1 \cdots y_d \partial_k f_{i} \bigg] \nonumber \\
=& \sum_{k=1}^d \sum_{\ell=1}^N 2 \mathsf{A}_{\xi} y_k (\partial_\ell f_k)^2 \partial_k f_i 
+ \sum_{k=1}^d y_k^{-1} \mathsf{A} \partial_k f_{i} + 2y_i^{-1} \mathsf{A} \partial_i f_{i} \\
&+ \sum_{k=1}^N \bigg[\sum_{j=1}^d \sum_{\ell=1}^N 2 \mathsf{A}_{\xi} y_j^2 \partial^2_{\ell k} f_{j} \partial_\ell f_{j} \partial_k f_{i} 
+ y_k [2\mathsf{A}_{s} + |x|^{-1}\mathsf{A}_{r}]\partial_k f_{i} + \mathsf{A} \partial_k^2 f_{i} \bigg]. \nonumber
\end{align}
This then relates to the operator $\msL\uab$ by way of the following result.

\begin{proposition} \label{connection systems}
Suppose $u$ is a whirl associated with the matrix field $\Q=\Q[\mathsf{f}]$ with $\mathsf{f}=(f_1, \dots, f_d)$ of class $\mathscr{C}^2$
$[$see \eqref{Qeven}-\eqref{Qodd}$]$. Then 
\begin{align}\label{Lwithdiv}
\msL[u;\A,\B] =& \, \nabla \A - \A \sum_{i=1}^d |\nabla f_i|^2 w^i 
+ \sum_{i=1}^d \frac{1}{ \mathscr{J} y_i^2} 
\dive[\mathscr{A}_i (y, \nabla \mathsf{f}) \nabla f_i] [w^i]^\perp \nonumber \\
&+ \sumell \sum_{i=1}^d  \frac{\partial_\ell f_i}{\mathscr{J}y_\ell} \dive[\mathscr{A}_i(y, \nabla \mathsf{f}) \nabla f_i] w^\ell + \B x.
\end{align}
The arguments of $\A, \B$ are $(r,s,\xi) = (||y||, ||y||^2, n+ \sum_{j=1}^d y_j^2|\nabla f_j|^2)$ and the coefficients 
$\mathscr{A}_i=\mathscr{A}_i(y, \nabla \mathsf{f})$ $($with $1 \le i \le d$$)$ are as in \eqref{nonlinearity A def}.
\end{proposition}

\begin{proof}
Starting from \eqref{LQ} and making use of \eqref{QtQl}, \eqref{QtQl-next}  we can rewrite $\msL[u;\A,\B]$ in terms of the components of the 
vector-map $\mathsf{f}$ as 
\begin{align}\label{Lwf}
&\msL[u; \A, \B] \, = \, \nabla \A + \B x + 2\A_\xi \times \nonumber \\
& \times \left\{ \sum_{k=1}^N \left[ \sum_{\ell=1}^N \left( y_k (\partial_\ell f_k)^2 
+ \sum_{i=1}^d y_i^2 \partial^2_{\ell k} f_i \partial_\ell f_i \right) \right] 
\left[\sum_{j=1}^d \left( \partial_k f_j [w^j]^\perp 
+ \sum_{i=1}^N y_j^2 \partial_i f_j \partial_k f_j \frac{w^i}{y_i} \right) \right] \right\} \nonumber \\
&+ \left[2\A_s + |x|^{-1}\A_r\right] \sum_{\ell=1}^N \left\{ \sum_{i=1}^d y_\ell \partial_\ell f_i [w^i]^\perp 
+ \sum_{k=1}^N \sum_{j=1}^d y_j^2 y_k \partial_\ell f_j \partial_k f_j \frac{w^\ell}{y_\ell}\bigg] \right\} \nonumber \\
&+ \A \left\{\sum_{\ell=1}^N \sum_{i=1}^d \left[  \left(\partial_\ell^2f_i + \Delta y_\ell \partial_\ell f_i\right)[w^i]^\perp 
- (\partial_\ell f_i)^2 w^i \right]+ \sum_{\ell=1}^d\frac{2}{y_\ell}\partial_\ell f_\ell [w^\ell]^\perp \right. \nonumber \\
&+ \left. \sumell \left[ \sum_{k=1}^N  \sum_{j=1}^d y_j^2\partial_\ell f_j \left( \partial^2_k f_j + \Delta y_k \partial_k f_j\right) 
+ \sum_{k=1}^d 2y_k \partial_\ell f_k \partial_k f_k \right] \frac{w^\ell}{y_\ell} \right\} .
\end{align}

Now referring to the expansion of the divergence operator prior to the proposition [{\it see} \eqref{divexp}] after a rearrangement of terms and a tedious but routine
set of calculations we arrive at the required conclusion. 
\end{proof}

The above result leads to two main consequences. The first underlines the role of the unconstrained system \eqref{ELwhirl} 
in relation to the solvability of the original system \eqref{EL1} and the second describes a stark simplification of the vector 
field $\msL[u;\A,\B]$ given that $u$ satisfies the restricted system \eqref{ELwhirl}.

\begin{proposition} \label{whirl full vs restrcited} If a whirl $u$ associated with the matrix field $\Q=\Q[\mathsf{f}]$ of class $\mathscr{C}^2$ is a solution 
to \eqref{EL1} then the vector-map $\mathsf{f}=(f_1, \dots, f_d)$ is a solution to \eqref{ELwhirl}. 
\end{proposition}

\begin{proof} Fixing $1 \le j \le d$ and taking the inner product of $\mathscr{L}[u; \A, \B]$ with $[w^j]^\perp$ by using the formulation in \eqref{Lwithdiv} 
and utilising the various orthogonality relations it is seen that  
\begin{align} \label{getting res PDE from full PDE}
\la \msL[u; \A, \B], [w^j]^\perp \ra =& \, \la \nabla \A, [w^j]^\perp \ra + \la \B x, [w^j]^\perp \ra - \la \A \sum_{i=1}^d |\nabla f_i|^2 w^i, [w^j]^\perp \ra \nonumber \\
&+ \frac{1}{\mathscr{J}} \sum_{i=1}^d \frac{1}{y_i^2} \la \dive[\mathscr{A}_i (y, \nabla \mathsf{f}) \nabla f_i] [w^i]^\perp, [w^j]^\perp \ra \nonumber \\
&+ \sum_{i=1}^d \sumell \la \frac{\partial_\ell f_i}{\mathscr{J}y_\ell} \dive[\mathscr{A}_i(y, \nabla \mathsf{f}) \nabla f_i] w^\ell, [w^j]^\perp \ra \nonumber \\
=& \,  \la \nabla \A, [w^j]^\perp \ra + \frac{1}{\mathscr{J}} \dive[\mathscr{A}_j (y, \nabla \mathsf{f}) \nabla f_j].
\end{align}
Now since $\A$ here is a function of $y=(y_1, \dots, y_N)$ an easy differentiation shows that its gradient $\nabla \A$
is a linear combination of the vectors $w^1, \dots, w^N$ and so $\la \nabla \A, [w^j]^\perp \ra \equiv 0$. As a 
result \eqref{getting res PDE from full PDE} simplifies further to $\la \msL[u; \A, \B], [w^j]^\perp \ra = 1/\mathscr{J} \dive[\mathscr{A}_j (y, \nabla \mathsf{f}) \nabla f_j]$.

Let $x=(x_1, \dots, x_n)$ be associated with the $2$-plane radial variables $y=(y_1, \dots, y_N)$. 
Consider the circle of radius $y_j$ given by $\gamma(t) = w^1 + \dots + w^j(t) + \dots + w^N$ ($0 \le t \le 2\pi$). Here $w^j(t) = y_j (0, \dots, 0, \cos t, \sin t, 0, \dots, 0)$ and except 
for $w^j$ all the other coordinates $w^\ell$ are independent of $t$. Then firstly all the points $x=\gamma(t)$ ($0 \le t \le 2\pi$) are associated with the same $y$ and secondly at the point 
$x=\gamma(t)$ we have $\dot \gamma(t) = d \gamma/dt(t) = [w^j(t)]^\perp$. Therefore by \eqref{getting res PDE from full PDE} and the PDE we have 
$\la \nabla \mcP, [w^j]^\perp \ra = 1/\mathscr{J} \dive[\mathscr{A}_j (y, \nabla \mathsf{f}) \nabla f_j]$ and hence substituting $x=\gamma(t)$, noting that the right-hand side is independent 
of $t$ and integrating over $0 \le t \le 2\pi$ we arrive at  
\begin{align}
\frac{2\pi}{\mathscr{J}(y)} \dive[\mathscr{A}_j (y, \nabla \mathsf{f}) \nabla f_j] &= \int_0^{2\pi} \la \nabla \mcP (\gamma(t)), [w^j(t)]^\perp \ra \, dt \\
&= \int_0^{2\pi} \la \nabla \mcP (\gamma(t)), \dot \gamma(t) \ra \, dt  = \int_0^{2\pi} \frac{d}{dt} \mcP(\gamma(t)) \, dt =0,  \nonumber
\end{align}
as required where the last identity follows from the closedness of $\gamma$. The proof is thus complete. 
\end{proof}

\begin{proposition}\label{reduction}
Assume $\mathsf{f} = (f_1,\dots,f_d)$ of class $\mathscr{C}^2$ is a solution to the system \eqref{ELwhirl}. Then denoting by $u$ the whirl 
associated with the matrix field $\Q=\Q[\mathsf{f}]$ we have 
\begin{equation}\label{Lsimp}
\msL[u;\A,\B] = \nabla\A  + \B x - \sum_{i=1}^d \A|\nabla f_i|^2 w^i .
\end{equation}
\end{proposition}

\begin{proof}
This follows from \eqref{Lwithdiv} upon noting that $\dive[\mathscr{A}_i(y, \nabla \mathsf{f}) \nabla f_i]=0$ by assumption for all $1 \le i \le d$.
\end{proof}



Proceeding forward recall that the overarching goal is to resolve the PDE $\msL\uab = \nabla\mcP$. Towards this end we consider the 
following general result before scrutinising the curl of the vector field in \eqref{Lsimp}.

\begin{lemma}\label{curlresult}
Consider the vector field $U(x) = \sum_{k=1}^N \bA_k(y) w^k$ for $\bA_k \in \mathscr{C}^1(\An,\R^{n\times n})$ for each $1 \leq k \leq N$. 
Then writing ${\bf W}_k = \nabla w^k = [\partial w_i^k/\partial x_j : 1 \le i, j \le n]$ we have 
\begin{equation}
{\rm curl}\, U = \sum_{k=1}^N \sumell \frac{1}{y_\ell} \left[\partial_\ell \bA_k w^k \otimes w^\ell - w^\ell \otimes \partial_\ell \bA_k w^k \right]
+ \sum_{k=1}^N (\bA_k {\bf W}_k - {\bf W}^t_k \bA_k^t).
\end{equation}
\end{lemma}

\begin{proof}
By linearity it suffices to consider only the case $U = {\bf A}_k w^k$. The conclusion then follows by a summation over $1 \le k \le N$. 
Towards this end by directly evaluating the curl and employing the product rule we have  
\begin{align*}
[{\rm curl}\,\bA_kw^k]_{ij} =& \, [\bA_kw^k]_{i,j} - [\bA_kw^k]_{j,i} \qquad \qquad \qquad \qquad \qquad \qquad 1 \leq i,j \leq n,  \\
=& \, \sumell \bigg[ [\partial_\ell \bA_kw^k ]_i\frac{w^\ell_j}{y_\ell} 
+ [\bA_k \nabla w^k]_{ij} - [\partial_\ell \bA_kw^k]_j\frac{w^\ell_i}{y_\ell} - [\bA_k \nabla w^k]_{ji} \bigg].
\end{align*}
This upon shifting to tensor notation immediately leads to the desired conclusion. 
\end{proof}

\begin{remark}{\em 
An easy inspection shows that ${\bf W}_k$ is the symmetric block diagonal matrix ${\bf W}_k=diag(0, \dots, 0, {\bf I}_2, 0, \dots, 0)$ 
with ${\bf I}_2$ as the $k^{th}$ block except for when $n=2d+1$ and $k=N$ in which case $\nabla w^N=diag(0, \dots, 0, 1)$.
}
\end{remark}

\begin{remark} \label{scalar curl remark} {\em 
In the case ${\bf A}_k(y)=\Gamma_k(y) {\bf I}_n$ with $U=\sum_{k=1}^N \Gamma_k(y) w^k$ and $\Gamma_k=\Gamma_k(y)$ 
suitable scalar functions, by using Lemma $\ref{curlresult}$ and $\bA_k {\bf W}_k - {\bf W}^t_k \bA_k^t=0$, we have
\begin{align}\label{easierVcurl}
{\rm curl}\,U =& \, \sum_{k=1}^N \sumell \frac{\partial_\ell \Gamma_k}{y_\ell} \left[w^k \otimes w^\ell - w^\ell \otimes w^k\right] \nonumber \\
=& \, \sum_{1 \le k < \ell \le N} \left( \frac{\partial_\ell \Gamma_k}{y_\ell} - \frac{\partial_k \Gamma_\ell}{y_k} \right) 
\left[w^k \otimes w^\ell - w^\ell \otimes w^k\right].
\end{align}
By virtue of the independence of the skew-symmetric tensors $[w^k \otimes w^\ell - w^\ell \otimes w^k]$ it follows by a continuity argument 
that ${\rm curl}\,U \equiv 0 \iff  \partial_\ell \Gamma_k/y_\ell - \partial_k \Gamma_\ell/y_k \equiv 0$ for all $1 \le k < \ell \le N$.}
\end{remark}

Returning now to \eqref{Lsimp} and by subtracting the gradient term $\nabla\A$ from both sides (and for the sake of uniformity 
in notation, extending the vector-map $\mathsf{f}$ in the case $n=2d+1$ to an $N$-vector by setting $f_N \equiv 0$) we consider the vector 
field $U = \msL\uab - \nabla\A$. 
This corresponds to the case in Remark \ref{scalar curl remark} with $\Gamma_k(y) = \B - \A|\nabla f_k|^2$ ($1 \leq k \leq N$) clearly of class $\mathscr{C}^1$. 
The ongoing analysis leading to \eqref{easierVcurl} then gives 
\begin{equation*}
{\rm curl}\, U = \sum_{1 \leq k < \ell \leq N} 
\bigg( \frac{\partial_\ell \B}{y_\ell} - \frac{\partial_k \B}{y_k} + \frac{\partial_k [\A|\nabla f_\ell|^2]}{y_k} - \frac{\partial_\ell [\A|\nabla f_k|^2]}{y_\ell} \bigg) 
\left[w^k \otimes w^\ell - w^\ell \otimes w^k \right] .
\end{equation*}

\begin{corollary}\label{cfcond}
The $\mathscr{C}^1$ vector field $U=\msL\uab - \nabla\A$ satisfies ${\rm curl}\, U \equiv 0$ iff 
\begin{equation}\label{curlrec}
y_k \partial_\ell \B - y_\ell \partial_k \B + y_\ell \partial_k [\A|\nabla f_\ell|^2] - y_k \partial_\ell [\A|\nabla f_k|^2] \equiv 0, 
\end{equation}
for all $1 \leq k < \ell \leq N$.
\end{corollary}


\section{Full Resolution of the System with $\A = H(r,s)$ and $\B=\B(r,s,\xi)$: The role of the Discriminant $\Delta$}
\label{Non variational WD Sec}

In this section we consider the nonlinear system $\Sigma[(u, \mathscr{P}), \A, \B]$ in \eqref{EL1} where 
we take $\A= H(r,s)$ where $H>0$ is of class $\mathscr{C}^2$ and $\B= \B(r,s,\xi)$ exactly as before. In 
this setting the operator \eqref{Lintro} and the associated PDE take the form 
\begin{align}\label{ELh}
\msL[u; H, \B] =& [\nabla u]^t \left\{ \dive[H(|x|, |u|^2) \nabla u] +\B(|x|, |u|^2, |\nabla u|^2) u \right\} \nonumber \\
=& [\nabla u]^t [\nabla u] \nabla H(|x|, |u|^2) + H(|x|, |u|^2) [\nabla u]^t \Delta u \nonumber \\
&+ \B(|x|, |u|^2, |\nabla u|^2) [\nabla u]^t u = \nabla \mcP. 
\end{align}
In line with the preceding analysis we also consider the unconstrained system \eqref{ELwhirl}, for the vector function 
$\mathsf{f}=(f_1, \dots, f_d)$, that in this context takes the form
\begin{align}\label{restrictDir}
\left\{
\begin{array}{ll}
\dive \left[ \mathscr{A}_\ell (y) \nabla f_\ell\right] = 0 \qquad &\text{ in } \An , \\
f_\ell \equiv g_\ell \qquad &\text{ on } (\partial\An)_\mathsf{D}, \qquad 1 \le \ell \le d,  \\
\mathscr{A}_\ell (y) \partial_\nu f_\ell = 0 \qquad &\text{ on } (\partial\An)_\mathsf{N}. 
\end{array}
\right. 
\end{align}
Here $\mathscr{A}_\ell(y) = y_\ell^2 H(z, z^2) \mathscr{J}(y)$, $g_\ell = 0$ at $z=a$ and $g_\ell = 2m_\ell\pi$ at $z=b$ with $m_\ell \in \mathbb{Z}$. 
Note that since $\A(r,s,\xi) = H(r,s)$ has no explicit $\xi$-dependence, unlike the original system \eqref{ELh}, 
here, the unconstrained system decouples and the $\ell^{th}$ PDE depends solely on the component $f_\ell$ rather 
than the full vector-map $\mathsf{f} = (f_1,\dots, f_d)$. This allows us to explicitly solve \eqref{ELh} in all dimensions  
which then leads to interesting consequences. 

\begin{theorem}\label{restrictedDirsolns}
Given $\mathsf{m} \in \mathbb{Z}^d$ the system \eqref{restrictDir} has the unique solution $\mathsf{f}=(f_1, \dots, f_d)$ given by 
\begin{equation}
\label{ODEDirsol}
\mathsf{f} = \mathsf{f}(y, \mathsf{m}) = 2 \mathsf{m}\pi \frac{\mathsf{H}(\|y\|)}{\mathsf{H}(b)} , \qquad \mathsf{H}(r) = \int_a^r \frac{dz}{z^{n+1}H(z,z^2)} .
\end{equation} 
\end{theorem}

\begin{proof}
The Dirichlet boundary condition $\mathsf{f}=\mathsf{g}$ on $(\partial\An)_\mathsf{D}$ is easily seen to be satisfied as a result of the normalisation 
and end-point conditions on $\mathsf{H}$. The Neumann boundary conditions $\mathscr{A}_\ell (y) \partial_\nu f_\ell = 0$ on $(\partial\An)_\mathsf{N}$ 
follow suit as a result of the quantity $\mathscr{J} \equiv 0$ on $(\partial\An)_\mathsf{N}$. Now referring to \eqref{ODEDirsol} a direct verification gives 
with $z=||y||$, 
\begin{equation}
\nabla \mathsf{f} = \left[ \frac{\partial f_i}{\partial y_j} = 2m_i\pi\frac{\dot{\mathsf{H}}(z)}{\mathsf{H}(b)}\frac{y_j}{z} : 1 \le i \le d, 1 \le j \le N \right] 
= \frac{2 \pi}{\mathsf{H}(b)} \frac{\mathsf{m} \otimes y}{z^{n+2} H(z,z^2)} .
\end{equation}
In even dimensions with $n =2d$ and $N = d$, we see by a direct calculation that,
\begin{align} \label{evenwhirl}
\dive \left[ \mathscr{A}_\ell(y) \nabla f_\ell \right] 
&= \sum_{j=1}^N \frac{\partial}{\partial y_j} \wc\left[\frac{y_j y_\ell^2}{||y||^{n+2}} \mathscr{J}(y)\right] \nonumber \\
&= \frac{\mathscr{J}(y)}{||y||^{n+2}}\wc y_\ell\left[dy_\ell - (2d+2)y_\ell + 2y_\ell +dy_\ell \right] =0.
\end{align}
In odd dimensions with $n = 2d+1$ and $N = d+1$, proceeding similarly and separating the first $y_1, \dots, y_d$ variables from the last variable 
$y_N$ in calculating the divergence, we see by a straightforward differentiation that,  
\begin{align} \label{oddwhirlk}
\dive \left[ \mathscr{A}_\ell(y) \nabla f_\ell \right] 
=& \sum_{j=1}^{d} \wc \frac{\partial}{\partial y_j} \bigg[\frac{y_j y_\ell^2}{||y||^{n+2}}\mathscr{J}(y)\bigg] 
+  \wc\frac{\partial}{\partial y_N}\left[\frac{y_Ny_\ell^2}{||y||^{n+2}}\mathscr{J}(y)\right] \nonumber \\
=& \frac{\mathscr{J}(y)}{||y||^{n+2}} \wc y_\ell\bigg[dy_\ell 
-\frac{(2d+3)}{||y||^2} y_\ell\sum_{i=1}^dy_i^2+ 2y_\ell + dy_\ell\bigg] \nonumber \\
& + \frac{\mathscr{J}(y)}{||y||^{n+2}}\wc y_\ell\left[y_\ell - \frac{(2d+3)}{||y||^2}y_\ell y_N^2\right] =0.
\end{align}
Having verified the solution to the PDE in \eqref{restrictDir} in both even and odd dimensions the assertion is justified and the proof is thus complete.
\end{proof}

Henceforth we shall write $\mcH(r) = \mathsf{H}(r)/\mathsf{H}(b)$. It is then evident that $\mcH$ is a solution to the linear 
ODE $d/dr [r^{n+1} H(r,r^2)\dot{\mcH}]=0$ on $a < r < b$. We write $\Q=\Q[f](y, \mathsf{m})$ for the matrix field 
associated with the solution $\mathsf{f}=\mathsf{f}(y, \mathsf{m})$ from Theorem \ref{restrictedDirsolns} 
[{\it see} \eqref{Qeven} and \eqref{Qodd}]. Thus it is plain that 
${\bf Q}[f](y, \mathsf{m}) = diag({\bf R}[ 2m_1 \pi \mcH(||y||)], \dots, {\bf R}[2m_d \pi \mcH(||y||)])$ when $n = 2d$ and 
${\bf Q}[f](y, \mathsf{m}) = diag({\bf R}[ 2m_1 \pi \mcH(||y||)], \dots, {\bf R}[2m_d \pi \mcH(||y||)], 1)$ when $n = 2d+1$. 
With these assumptions in place the action of the differential operator $\msL$ on the map $u=\Q[\mathsf{f}](y) x$ 
after the subtraction of $\nabla H$ [{\it cf.} \eqref{Lsimp}] can be written   
\begin{align}\label{LDir}
U &= \msL[u; H,\B] - \nabla H = - H \sum_{i=1}^d |\nabla f_i|^2 w^i + \B \left(z, z^2, n+\sum_{j=1}^d y_j^2 |\nabla f_j|^2\right) x.
\end{align}

We now turn to the curl of the vector field \eqref{LDir} in anticipation of solving the PDE $\msL[u; H,\B] =\nabla\mcP$. 
Here we use Lemma \ref{curlresult} and Remark \ref{scalar curl remark} with the choice of functions 
\begin{equation}
\Gamma_k(y) = - 4m_k^2\pi^2 \dot{\mcH}^2 H + \B \left( z, z^2, n + \sum_{j=1}^d 4 m_j^2 \pi^2 y_j^2 \dot{\mcH}^2 \right), 
\qquad 1 \le k \le N,  
\end{equation}
noting that $|\nabla f_j|^2=4m_j^2 \pi^2 \dot{\mcH}^2$. Remark \ref{scalar curl remark} and Corollary \ref{cfcond} direct us to compute the expressions 
$\partial_\ell \Gamma_k / y_\ell - \partial_k \Gamma_\ell / y_k$. Towards this end we first observe that 
\begin{align}
\partial_\ell \Gamma_k(y) 
=& \, 8\pi^2 \B_\xi \bigg[\dot{\mcH} \ddot{\mcH} \frac{y_\ell}{||y||} \sum_{j=1}^d m_j^2 y_j^2 
+ \dot{\mcH}^2 m_\ell^2 y_\ell \bigg] + \bigg( 2||y|| \B_s + \B_r \bigg) \frac{y_\ell}{||y||} \nonumber \\
&- 4m_k^2\pi^2 \left[ 2\dot{\mcH} \ddot{\mcH} H + \dot{\mcH}^2 \frac{dH}{dr} \right] \frac{y_\ell}{||y||}, 
\end{align}
where we have abbreviated the arguments of $\B_r, \B_s$ and $\B_\xi$. Therefore it follows that 
\begin{align}\label{curlhcond1}
\frac{\partial_\ell \Gamma_k}{y_\ell} - \frac{\partial_k \Gamma_\ell}{y_k} =& 
- \frac{1}{y_k} \bigg\{ 8\pi^2 \B_\xi \bigg[\dot{\mcH} \ddot{\mcH} \frac{y_k}{||y||} \sum_{j=1}^d m_j^2 y_j^2 
+ \dot{\mcH}^2 m_k^2 y_k \bigg] + \bigg( 2||y|| \B_s + \B_r \bigg) \frac{y_k}{||y||} \nonumber \\
& \qquad - 4m_\ell^2\pi^2 \left[ 2\dot{\mcH} \ddot{\mcH} H + \dot{\mcH}^2 \frac{dH}{dr} \right] \frac{y_k}{||y||} \bigg\} 
\nonumber \\
&+ \frac{1}{y_\ell} \bigg\{ 8\pi^2 \B_\xi \bigg[\dot{\mcH} \ddot{\mcH} \frac{y_\ell}{||y||} \sum_{j=1}^d m_j^2 y_j^2 
+ \dot{\mcH}^2 m_\ell^2 y_\ell \bigg] + \bigg( 2||y|| \B_s + \B_r \bigg) \frac{y_\ell}{||y||} \nonumber \\
& \qquad - 4m_k^2\pi^2 \left[ 2\dot{\mcH} \ddot{\mcH} H + \dot{\mcH}^2 \frac{dH}{dr} \right] \frac{y_\ell}{||y||} \bigg\} 
\nonumber \\
=& \, 4(m_\ell^2 - m_k^2) \pi^2 \left[ 2\dot{\mcH}^2 \B_\xi + \frac{1}{r} \dot{\mcH}^2 \frac{dH}{dr} + \frac{2}{r} \dot{\mcH} \ddot{\mcH} H \right] 
\nonumber \\ 
=& \, 4(m_k^2 - m_\ell^2) \pi^2 \frac{\dot{\mcH}^2}{||y||^2} \Delta(H,\B). 
\end{align}
Here we have introduced the {\it discriminant} $\Delta(H,\B) := 2(n+1)H + rH_r + 2r^2[H_s - \B_\xi]$ associated with $(H, \B)$ and 
the putative map $u$. Now \eqref{easierVcurl} and Corollary~\ref{cfcond} give
\begin{equation} \label{curlhdisc}
{\rm curl}\,[\msL (u; H, \B) - \nabla H] = \frac{4\pi^2\dot{\mcH}^2}{||y||^2} \Delta(H, \B) 
\sum_{1\leq k < \ell \leq N}\hspace{-2mm} \left(m_k^2 - m_\ell^2\right) \left[w^k\otimes w^\ell - w^\ell \otimes w^k \right] .
\end{equation}
The following theorem gives a complete characterisation of all whirl solutions to the system \eqref{EL1}-\eqref{Lintro} 
[with $\mathscr{L}$ as in \eqref{ELh}] pointing at an interesting dimensional parity.

\begin{theorem}\label{Dirsoln}
A whirl $u$ associated with the matrix field $\Q=\Q[\mathsf{f}] \in \mathscr{C}^2(\overline{\mathbb{A}}_n,\SO(n))$ and satisfying the boundary condition 
$\Q \equiv {\bf I}$ on $(\partial \An)_\mathsf{D}$ is a solution to the system $\Sigma[(u, \mathscr{P}); H, \B]$ in \eqref{EL1} with $\msL[u; H, \B]$ as 
in \eqref{ELh} if and only if the following hold.
\begin{enumerate}
\item[$\bullet$] If $\Delta(H, \B) \not\equiv 0$ then depending on the dimension being even or odd, we have:
\begin{enumerate}
\item[$(i)$] $n=2d$: $\Q=\Q[\mathsf{f}](y)= diag({\bf R}[2m_1\pi\mcH(||y||)], \dots , {\bf R}[2m_d\pi\mcH(||y||)])$ 
with $\mathsf{m}=(m_1, \dots, m_d) \in \mathbb{Z}^d$ satisfying $|m_1| = \dots = |m_d|$.
\item[$(ii)$] $n=2d+1$: $\Q=\Q[\mathsf{f}] \equiv \In$ corresponding to $m_1 = \dots = m_d = 0$.
\end{enumerate}
\item[$\bullet$] If $\Delta(H, \B) \equiv 0$ then depending on the dimension being even or odd, we have:
\begin{enumerate}
\item[$(i)$] $n=2d$: $\Q=\Q[\mathsf{f}] = diag({\bf R}[2m_1\pi\mcH(||y||)], \dots , {\bf R}[2m_d\pi\mcH(||y||)])$. 
\item[$(ii)$] $n=2d+1$: $\Q=\Q[\mathsf{f}] = diag({\bf R}[2m_1\pi\mcH(||y||)], \dots , {\bf R}[2m_d\pi\mcH(||y||)],1)$. 
\end{enumerate}
In either case $\mathsf{m}=(m_1, \dots, m_d) \in \mathbb{Z}^d$ and there is no further restrictions needed on $m_1,\dots, m_d$.
\end{enumerate}
\end{theorem}

\begin{proof} We shall split the proof into two parts justifying the necessity and the sufficiency arguments separately.  \\
({\it Necessity.}) By Proposition \ref{whirl full vs restrcited} if a whirl associated with the matrix field $\Q=\Q[\mathsf{f}]$ is a solution to 
the system $\Sigma[(u, \mathscr{P}); H, \B]$ then the vector-map $\mathsf{f}$ must be a solution to \eqref{ELwhirl}, or more specifically, 
here, to \eqref{restrictDir}. Therefore $\mathsf{f}$ must be exactly as described by Theorem \ref{restrictedDirsolns}. It suffices now to use 
Proposition \ref{reduction}, the curl analysis in Section \ref{Restricted and Curl Section} and the calculation leading to \eqref{curlhdisc} to 
get a complete characterisation of $\mathsf{f}$ and to do so we proceed by considering the cases $\Delta(H, \B) \not\equiv 0$ 
and $\Delta(H, \B) \equiv 0$ separately.

If $\Delta(H, \B) \not\equiv 0$ then by virtue of the independence of the tensors $w^k\otimes w^\ell - w^\ell \otimes w^k$ 
($1 \le \ell < k \le N$) and the fact these tenors vanish at most on coordinate hyperplanes, \eqref{curlhdisc} and a basic 
continuity argument gives ${\rm curl}\,U \equiv 0 \iff m_1^2 = \dots = m_N^2$. 
In the case $n=2d$ this gives the conclusion in $(i)$ and in the case $n=2d+1$ this gives $m_1 = \cdots = m_N = 0$ and 
so $\Q \equiv \In$ as stated in $(ii)$. If $\Delta(h,\B) \equiv 0$ then again by \eqref{curlhdisc} ${\rm curl}\,U \equiv 0$ 
irrespective of the choice see there can be no restriction on the integers $m_\ell$. \\
({\it Sufficiency.}) We shall do this only for the case $\Delta(H,\B) \equiv 0$ as the case $\Delta(H,\B) \not\equiv 0$ is straightforward. 
Towards this end we assume hereafter that $\Delta(H,\B) \equiv 0$ and show that $\mathscr{L}[u; H, \B] = \nabla \mathscr{P}$. We 
claim that $U=\msL[u; H, \B]-\nabla H = \nabla \mathscr{R}(|x|,|\bH x|^2)$ for a suitable choice of $\mathscr{R}=\mathscr{R}(r, z)$ 
of class $\mathscr{C}^2$. Here $\bH$ stands for the $n \times n$ skew-symmetric matrix $\bH = diag(2m_1\pi\bJ,\dots,2m_d\pi\bJ)$ 
or $\bH = diag(2m_1\pi\bJ,\dots,2m_d\pi\bJ,0)$ depending as to whether $n=2d$ or $n = 2d+1$. Indeed assuming the claim to be 
true a direct calculation and comparison with $U$ leads to  
\begin{align} \label{candidateP}
\nabla\mathscr{R}(|x|, |\bH x|^2) &= \mathscr{R}_r(r, |\bH x|^2) \Theta - 2r\mathscr{R}_z(r, |\bH x|^2)\bH^2 \Theta \nonumber \\
&= r \B(r, r^2, n + \dot \mcH^2 |\bH x|^2) \Theta + r H \dot\mcH^2 \bH^2 \Theta = U, 
\end{align}
where $\mathscr{R}_r$, $\mathscr{R}_z$ denote the derivatives of $\mathscr{R}$ in the first and second arguments respectively. 
Thus the second equality in \eqref{candidateP}  would be valid [{\it cf.} \eqref{LDir}] provided that $\mathscr{R}_r(r, z) = r\B(r, r^2, n+\dot\mcH^2z)$ and 
$\mathscr{R}_z(r,z) = -H(r,r^2)\dH^2/2$. Let us thus turn on to constructing $\mathscr{R}$. To this end let $\mathscr{B}=\{(r, z) : \mbox{$r=|x|$, $z=|\bH x|^2$ 
with $x \in \Xn$}\}$. Then $\mathscr{B} \subset [a, b] \times \mathbb{R}$ is seen to be simply-connected; as a matter of fact, denoting by 
$\underline{m}, \overline{m} \ge 0$ the minimum and maximum eigenvalues of the diagonal matrix $\bH^t \bH$ 
respectively it is easily seen that $\mathscr{B}=\{(r, z) : a<r<b, 0 \le \underline{m} r^2  \le z \le \overline{m} r^2\}$. 
Since $\Delta (H, \B) \equiv 0$ it is not difficult to see that 
$\partial_z \mathscr{R}_r (r, z) - \partial_r \mathscr{R}_z (r, z) = \partial_z [r\B(r, r^2, n+\dot\mcH^2z)] + \partial_r [H(r,r^2)\dH^2/2] \equiv 0$ in $\mathscr{B}$. 
As a result the $1$-form $\omega = r\B(r, r^2, n+\dot\mcH^2z) \, dr - H(r,r^2)\dH^2/2\, dz$ is closed in $\mathscr{B}$ and hence exact in view of $\mathscr{B}$ 
being simply-connected. Thus $\omega = d \mathscr{R}$ for a function (a $0$-form) $\mathscr{R}=\mathscr{R}(r,z)$ of class $\mathscr{C}^2$. To describe 
$\mathscr{R}$ more specifically pick a base point $(r^\star, z^\star)$ in $\mathscr{B}$ and let $\gamma$ be any piecewise continuously differentiable Jordan 
curve in $\mathscr{B}$ connecting $(r^\star, z^\star)$ to $(r,z)$ and set 
\begin{equation} \label{nablaU}
\mathscr{R}(r,z) = \int_\gamma \omega = \int_\gamma   r\B(r, r^2, n+\dot\mcH^2z) \, dr - H(r,r^2)\dH^2(r) /2\, dz, \qquad (r, z) \in \mathscr{B}. 
\end{equation}
The integral is seen to be independent of the choice of $\gamma$ and hence well-defined. The function $\mathscr{R}$ is of class $\mathscr{C}^2$ in the interior of 
$\mathscr{B}$ with continuously differentiable tangential gradients on the upper and lower boundary curves of $\mathscr{B}$. One can thus verify that 
\eqref{candidateP} holds (both for $(r,z)=(|x|, |\bH x|^2)$ in the interior of $\mathscr{B}$ and the upper and lower boundary curves).
Thus $U=\nabla\mathscr{R}(|x|, |\bH x|^2)$ and the proof is complete.
\end{proof}


\section{Infinitely many whirl solutions to $\Sigma[(u, \mcP); \A, \B]$ in even dimensions}
\label{Infinitude in even d Sec}

In this last section we prove the existence of an infinitude of solutions to the original system $\Sigma[(u, \mcP), \A, \B]$ in even dimensions. 
In terms of our earlier notation here $n=2d$ with $d=N$ and the Dirichlet boundary condition will be chosen 
$\mathsf{g}=2m \pi \mathsf{1} \chi_{z=b}$ with $m \in \mathbb{Z}$. Here $\chi_{z=b}$ is the characteristic function of the set $\{z=b\}$.

\begin{theorem}\label{whirlsoln}
Let $n = 2d$ and $\mathsf{m} = m \mathsf{1} \in \mathbb{Z}^d$. Then \eqref{ELwhirl} admits the unique solution 
$\mathsf{f}(y;\m) = (f_1,\dots,f_d)=\mcG(\|y\|; m) \mathsf{1}$ where $\mcG=\mcG(r; m) \in \mathscr{C}^2[a,b]$ 
is the unique solution to the two point boundary-value problem 
\footnote{The existence of such solutions $\mcG$ with the required $\mathscr{C}^2$-regularity can be established as in \cite{MsonTb}.}
\begin{align}\label{BVPG}
\begin{dcases}
\dfrac{d}{dr}\left[r^{n+1}\A(r, r^2, n+r^2\dot{\mcG}^2)\dot{\mcG}\right] = 0 , \quad a < r < b , \\
\mcG(a) = 0, \\
\mcG(b) = 2m\pi.
\end{dcases}
\end{align}
\end{theorem}

\begin{proof}
The boundary conditions on $(\partial\An)_\mathsf{D}$ in \eqref{ELwhirl} follow from the imposed end-point conditions on $\mcG$ in \eqref{BVPG}. 
Now in order to verify the PDE in \eqref{ELwhirl} we first observe that, 
\begin{equation} \label{df}
\nabla \mathsf{f} = \left[ \frac{\partial f_i}{\partial y_j} : 1 \le i \le d, 1 \le j \le N \right] 
= \dot{\mcG} \frac{{\mathsf{1} } \otimes y}{||y||} \ \implies \ 
\sum_{\ell = 1}^N y_\ell^2|\nabla f_\ell|^2 
= \sum_{\ell = 1}^N y_\ell^2 \dot{\mcG}^2 = ||y||^2\dot{\mcG}^2 .
\end{equation}
Now upon noting that if $n=2d$ we have $y=(y_1, \dots, y_d)$ and $\dive = \sum_{j=1}^d \partial/\partial y_j$, we can proceed directly and write  
\begin{align}\label{divformula}
\dive [\mathscr{A}_i (y, \nabla \mathsf{f}) \nabla f_i ] 
&= \sum_{j=1}^N \frac{\partial}{\partial y_j}\left[\mathscr{A}_i(y, \nabla \mathsf{f}) \frac{\partial f_i}{\partial y_j} \right] 
=\sum_{j=1}^d \frac{\partial}{\partial y_j}\left[y_i^2 \A \mathscr{J}(y)  \dot{\mcG} \frac{y_j}{||y||} \right] \nonumber \\
&= \frac{y_i^2 \mathscr{J}(y)}{||y||} \bigg\{ r\dot{\mcG}\frac{d}{dr}\A +r\A\ddot{\mcG} + (2d+1)\A\dot{\mcG} \bigg\} \nonumber \\
&=\frac{y_i^2 \mathscr{J}(y)}{||y||^{n+1}} \frac{d}{dr}\left[r^{n+1}\A(r,r^2,n+r^2\dot{\mcG}^2)\dot{\mcG}\right], 
\end{align}
where in the first two lines we have written $\A = \A(r,r^2,n+r^2\dot{\mcG}^2)$. It is now plain that if $\mcG$ is a solution to the ODE in 
\eqref{BVPG} then the vector $\mathsf{f}$ satisfies the above PDE. This therefore completes the proof.
\end{proof}

We turn now to the system \eqref{EL1} and prove the multiplicity result announced at the start of the section. Indeed we prove that 
for each $m \in \mathbb{Z}$ the whirl map $u=u(x; m)$ with $\rad(x) = |x|$, $\sph (x) = \Q[\mathsf{f}](y, \mathsf{m}) \Theta$ where 
$\Q=\Q[\mathsf{f}](y,\m) = \matrixexp\{\mcG(||y||;m)\bH\} = diag({\bf R}[\mcG(||y||; m)], \dots, {\bf R}[\mcG(||y||; m)])$ and $\mcG$ 
is as in Theorem \ref{whirlsoln} serves as a solution to \eqref{EL1}. Here we write $\bH$ for the $\mathfrak{so}(n)$ matrix 
$\bH = diag(\bJ,\dots,\bJ)$.

\begin{theorem} \label{multiplicity even F}
For $n \geq 2$ even and $m \in \mathbb{Z}$ let $u = r\matrixexp\{\mcG(r;m)\bH\}\Theta$ where $\mcG\in\mathscr{C}^2[a,b]$ 
is the solution to \eqref{BVPG} and $\bH = diag(\bold{J},\dots,\bold{J})$. Then $\msL\uab = \nabla\mcP$ where the pressure 
field takes the form $\mcP = \A + \mathsf{G}$. Here the radial scalar-valued function $\mathsf{G} = \mathsf{G}(r)$ is chosen 
such that $\nabla G = r[\B(r,r^2,n+r^2\dG^2) - \A(r,r^2,n+r^2\dG^2)\dG^2]\Theta$. As a result the system \eqref{EL1} has an 
infinitude of $\mathscr{C}^2$ solutions.
\end{theorem}

\begin{proof}
The boundary conditions in \eqref{EL1} follow immediately from those of $\mcG$ in \eqref{BVPG} and as seen earlier $\det\nabla u = 1$. 
It thus remains to prove that $\msL\uab$ is a gradient field. Taking $\mathsf{f}$ as in Theorem \ref{whirlsoln} it follows after an application 
of Proposition~\ref{reduction} that $\msL\uab = \nabla \A+r[\B(r,r^2,n+r^2\dG^2) - \A(r,r^2,n+r^2\dG^2)\dG^2]\Theta$. From this 
the description of $\mcP = \mcP(x)$ as in the statement of the theorem follows and the proof is complete.
\end{proof}


\section{Appendix}
\label{Derivation-EL}

In this appendix we give a short derivation of the equations of first variation for the total elastic energy integral 
\eqref{I-W-energy} subject to the incompressibility constraint (that is, the system \eqref{PDE1}-\eqref{PDE1-var}). 
The argument is known among the experts and is given here for the sake convenience of the reader. 
Towards this end we pick a map $u=(u_1, \dots, u_n)$. As is standard we derive the equations under the assumption 
of sufficient smoothness of $u$. Moreover $\det \nabla u \equiv1$ and we assume that $u$ is injective on $\Omega$. 
Setting $U=u(\Omega) \subset \mathbb{R}^n$ it follows from the invariance of domain that $U$ is open. We denote the inverse of $u$ by $u^{-1}$.  
\footnote{Note that the whirl solutions $u=\Q[\mathsf{f}](y) x$ constructed in the paper have both the required 
degree of smoonthness and are injective.}

Pick a smooth compactly supported vector field $v \in \mathscr{C}_0^\infty(U, \mathbb{R}^n)$ and assume that  
$\dive\,v = {\rm tr}\,\nabla v=0$ in $U$. Consider the integral curves of the vector field $v$ in $U$, i.e., for every 
$y \in U$ and $t \in \mathbb{R}$ let $\Upsilon=\Upsilon(y, t)$ denote the solution to the initial value problem 
\begin{align}\label{Flow}
\begin{dcases}
\dfrac{d}{dt} \Upsilon (y,t) = v(\Upsilon(y, t)), \qquad t \in \mathbb{R}, \\
\Upsilon(y, 0) = y.
\end{dcases}
\end{align}
Now recalling the relations $d(\det P)/d P = {\rm cof}\,P$ and $P  [{\rm cof}\,P]^t =  [{\rm cof}\,P]^t P= (\det \, P) {\bf I}_n$, a straightforward 
differentiation gives  
\begin{align}\label{detUpsilon}
\frac{d}{dt} \det \nabla_y \Upsilon &= \sum_{ij} [{\rm cof}\,\nabla_y\Upsilon]_{ij} \frac{d}{dt} [\nabla_y \Upsilon]_{ij} \nonumber \\
&= \sum_{ij} [{\rm cof}\,\nabla_y\Upsilon]_{ij} [\nabla_y v(\Upsilon)]_{ij}  
= \sum_{ijk} [{\rm cof}\,\nabla_y \Upsilon]_{ij} [\nabla_y v]_{ik} [\nabla_y \Upsilon]_{kj} \nonumber \\
&= [\nabla_y \Upsilon] [{\rm cof}\,\nabla\Upsilon]^t : [\nabla_y v]^t = (\det \nabla_y \Upsilon) \, \dive\,v. 
\end{align}
Thus in particular as the vector field $v$ is chosen to be divergence free we have 
\begin{equation}\label{detUpsilon-one}
\frac{d}{dt} \det \nabla_y \Upsilon = (\det \nabla_y \Upsilon) \, \dive\,v =0. 
\end{equation}
Next as by \eqref{Flow} we have $\det \nabla_y \Upsilon(y, 0) =1$ it follows from \eqref{detUpsilon-one} that $\det \nabla_y \Upsilon (y,t) =1$ 
for all $y \in U$, $t \in \mathbb{R}$. Let us now set $u_t(x) = \Upsilon(u(x), t)$ for $x$ in $\Omega$ and $t \in \mathbb{R}$. A basic calculation 
then gives 
\begin{equation*}
\det \{\nabla_x u_t (x)\} = \det \{\nabla_x [\Upsilon(u(x), t)] \} = \det \{[\nabla_y \Upsilon (u(x), t)] [\nabla_x u(x)] \} = 1.
\end{equation*}
Furthermore by an easy inspection $u_t=u$ near the boundary $\partial \Omega$ whilst $u_0=u$. As a result $u_t$ constitutes a one parameter 
family of incompressible deformations in $\mathscr{A}^p_\varphi (\Omega)$ passing through $u$ at $t=0$. Hence by referring to the energy 
integral \eqref{I-W-energy}, for $u$ to be an energy extremiser, we must have 
\begin{equation} \label{I-W-variation-energy}
\frac{d}{dt} \mathbb{E}[u_t] \bigg|_{t=0} = \frac{d}{dt} \int_\Omega W(x, u_t(x), \nabla u_t(x)) \, dx \bigg|_{t=0} =0.  
\end{equation}

A direct calculation and making use of \eqref{Flow} now gives 
\begin{align} \label{det-EL-calc}
& \frac{d}{dt} \mathbb{E}[u_t] \bigg|_{t=0} \nonumber \\
&= \int_\Omega \left( \left\langle W_u (x,u(x),\nabla u(x)), \frac{du_t}{dt} \right\rangle \bigg|_{t=0} 
+ W_\xi (x,u(x),\nabla u(x)) : \frac{d \nabla u_t}{dt} \bigg|_{t=0} \right) dx \nonumber \\
&= \int_\Omega \left( \langle W_u (x,u(x),\nabla u(x)), [v \circ u] (x) \rangle + W_\xi (x,u(x),\nabla u(x)) : \nabla [v \circ u](x) \right) dx. 
\end{align}
Let us set $W_u(x) = W_u(x, u(x), \nabla u(x))$ and $W_\xi(x)=W_\xi(x, u(x), \nabla u(x))$ for brevity. Then from \eqref{det-EL-calc} 
and after a change of variables (using the invertibility of $u$) we have  
\begin{align} \label{det-EL-calc-L}
\frac{d}{dt} \mathbb{E}[u_t] \bigg|_{t=0}  
&= \int_U \left( \langle W_u (u^{-1}(y)), v(y) \rangle + W_\xi (u^{-1}(y)) : [\nabla_y v](y) [\nabla_x u](u^{-1}(y))  \right) dy \nonumber \\
&= \int_U \left( \langle W_u (u^{-1}(y)), v(y) \rangle + W_\xi (u^{-1}(y)) [\nabla_x u]^t (u^{-1}(y)) : [\nabla_y v](y) \right) dy.
\end{align}
Next let us denote the integral on the right in \eqref{det-EL-calc-L} by $L(v)$, that is, let 
\begin{align} \label{L-def-U}
L(v) = \int_U \left( \langle W_u (u^{-1}(y)), v(y) \rangle + W_\xi (u^{-1}(y)) [\nabla_x u]^t (u^{-1}(y)) : [\nabla_y v] (y) \right) dy.
\end{align}
Then it is easily seen that there exists $c>0$ (depending on $u$, $W$ but independent of $v$) such that for all vector fields 
$v \in \mathscr{C}_0^1(U, \mathbb{R}^n)$ we have 
\begin{equation}
|L(v)| \le c \left[ ||v||_{L^\infty(U, \mathbb{R}^n)} + ||\nabla_y v||_{L^\infty(U, \mathbb{R}^{n \times n})} \right].
\end{equation} 
Thus $L$ is a bounded linear functional on 
$\mathscr{C}_0^1(U, \mathbb{R}^n)$. As from \eqref{I-W-variation-energy} and \eqref{det-EL-calc-L} we have $L(v) = 0$ 
for when $\dive\,v=0$ it then follows ({\it see} Section 1.4 and Proposition 1.1 in \cite{TeR}) that there exists $p \in \mathscr{D}'(U)$ 
such that $L = -\nabla p$. In particular we can write 
\begin{equation} \label{L-rep-grad}
L(v) = - (\nabla p, v ) = (p, \dive\,v) = (p, {\rm tr}\,[\nabla v]).
\end{equation} 
Now a reference to \eqref{L-def-U} and an application of the integration by parts formula on the second term in the integral together 
with \eqref{L-rep-grad} gives
\begin{align}
 \int_U \left( \langle W_u (u^{-1}(y)) - \dive\,\{W_\xi (u^{-1}(y)) [\nabla_x u]^t (u^{-1}(y))\}, v(y) \rangle  \right) dy = -(\nabla p, v).
\end{align}
This in particular implies that $\nabla p$ can be represented by an integrable function (in fact continuous) on $U$.   
Next let us take $\phi \in \mathscr{C}_0^\infty(\Omega; \mathbb{R}^n)$ and set $v=\phi \circ u^{-1}$. 
Then a straightforward differentiation results in
\begin{equation} \label{grad-v-calc}
\nabla_y v (y) = [\nabla_x \phi] (u^{-1}(y)) [\nabla_y u^{-1}](y) = [\nabla_x \phi] (u^{-1}(y)) [\nabla_x u]^{-1}(u^{-1}(y)).
\end{equation}
Hence substitution in the integral on the right \eqref{L-def-U} and changing variables by transforming back to $\Omega$ gives 
\begin{align} \label{L-W-form}
L(v)
&= \int_U \left( \langle W_u (u^{-1}(y)), \phi (u^{-1}(y)) \rangle + W_\xi (u^{-1}(y)) : [\nabla_x \phi] (u^{-1}(y)) \right) dy \nonumber \\
&= \int_\Omega \left( \langle W_u(x), \phi(x) \rangle + W_\xi(x) : [\nabla_x \phi](x) \right) dx.  
\end{align}
Likewise substitution in \eqref{L-rep-grad} and a similar argument after setting $\mcP=p \circ u$ results in
\begin{align} \label{L-q-form}
L(v) = - (\nabla p, v) 
&=\int_U p(y) \, \dive\,v(y) dy
= \int_U p(y)  \, {\rm tr} \{ [\nabla_x \phi] (u^{-1}(y) [\nabla_x u]^{-1} (u^{-1}(y)) \} \, dy \nonumber \\
&= \int_\Omega (p \circ u) (x) \, {\rm tr} \{ [\nabla_x \phi] [{\rm cof}\,\nabla_x u]^t \} \, dx \nonumber \\
&= \int_\Omega \mcP(x) \, [{\rm cof}\,\nabla_x u] : [\nabla_x \phi] \, dx.
\end{align}
Finally equating \eqref{L-W-form} and \eqref{L-q-form} leads to 
\begin{equation}
\int_\Omega \left( \langle W_u(x), \phi \rangle + (W_\xi(x) - \mcP(x) \, [{\rm cof}\,\nabla_x u] : [\nabla_x \phi] \right) dx = 0, 
\end{equation}
which after taking into account the arbitrariness of $\phi \in \mathscr{C}_0^\infty(\Omega, \mathbb{R}^n)$ and standard 
arguments formally results in \eqref{PDE1-var}.


\qquad \\
{\bf Acknowledgement.} AT acknowledges support from EPSRC grant EP/V027115/1 Topology of Sobolev Spaces and Quasiconvexity.

\end{document}